\documentclass[table,svgnames,english,oneside]{smfart}

\title[Formalization of non-Archimedean Functional Analysis]{Formalization of non-Archimedean Functional Analysis 1: Spherically Complete Spaces}
\usepackage[dvipsnames]{xcolor}
\usepackage{letltxmacro}
\LetLtxMacro\oldttfamily\ttfamily
\DeclareRobustCommand{\ttfamily}{\oldttfamily\csname ttsize\endcsname}
\newcommand{\setttsize}[1]{\def\ttsize{#1}}%
\usepackage{orcidlink,amsmath,mathtools,amssymb,enumitem,geometry,graphicx,float,tikz-cd,fontawesome5}
\usepackage[format=plain,justification=centering]{caption}
\usepackage{mathrsfs}

\DeclareSymbolFont{bbold}{U}{bbold}{m}{n}
\DeclareMathSymbol{\bbmu}{\mathord}{bbold}{"16}
\usepackage{xparse}
\ExplSyntaxOn
\NewDocumentCommand{\definealphabet}{mmmm}
 {% #1 = prefix, #2 = command, #3 = start, #4 = end
  \int_step_inline:nnn { `#3 } { `#4 }
   {
    \cs_new_protected:cpx { #1 \char_generate:nn { ##1 }{ 11 } }
     {
      \exp_not:N #2 { \char_generate:nn { ##1 } { 11 } }
     }
   }
 }
\ExplSyntaxOff
\definealphabet{bb}{\mathbb}{A}{Z}
\definealphabet{cal}{\mathcal}{A}{Z}
\definealphabet{frak}{\mathfrak}{A}{z}
\definealphabet{rm}{\mathrm}{A}{z}
\definealphabet{bf}{\mathbf}{A}{Z}
\definealphabet{scr}{\mathscr}{A}{Z}
\definealphabet{tt}{\mathtt}{A}{Z}
\let\opn\operatorname
\newcommand{\lto}{\longrightarrow}
\newcommand{\norm}[1]{\left\lVert#1\right\rVert}

\usepackage{hyperref}
\usepackage[nameinlink]{cleveref}
\Crefname{theorem}{Theorem}{Theorem}
\Crefname{conjecture}{Conjecture}{Conjectures}
\Crefname{lemma}{Lemma}{Lemmas}
\Crefname{definition}{Definition}{Definitions}
\Crefname{remark}{Remark}{Remarks}
\Crefname{proposition}{Proposition}{Propositions}
\Crefname{corollary}{Corollary}{Corollaries}
\Crefname{equation}{}{}
\Crefname{item}{}{}
\Crefname{example}{Example}{Examples}
\Crefname{claim}{Claim}{Claims}
\Crefname{assertion}{Assertion}{Assertions}
\Crefname{proof}{Proof}{Proofs}
\Crefname{condition}{Condition}{Conditions}
\Crefname{question}{Question}{Questions}
\newtheorem{theorem}{Theorem}[section]

\newtheorem{example}[theorem]{Example}
\AddToHook{env/example/begin}{\crefalias{theorem}{example}}
\newtheorem{lemma}[theorem]{Lemma}
\AddToHook{env/lemma/begin}{\crefalias{theorem}{lemma}}
\newtheorem{remark}[theorem]{Remark}
\AddToHook{env/remark/begin}{\crefalias{theorem}{remark}}
\newtheorem{proposition}[theorem]{Proposition}
\AddToHook{env/proposition/begin}{\crefalias{theorem}{proposition}}
\newtheorem{definition}[theorem]{Definition}
\AddToHook{env/definition/begin}{\crefalias{theorem}{definition}}

\AddToHook{env/conjecture/begin}{\crefalias{theorem}{conjecture}}
\newtheorem{corollary}[theorem]{Corollary}
\AddToHook{env/corollary/begin}{\crefalias{theorem}{corollary}}

\AddToHook{env/question/begin}{\crefalias{theorem}{question}}

\AddToHook{env/assumption/begin}{\crefalias{theorem}{assumption}}
\AddToHook{env/theorem/before}{\addvspace{3pt}}
\AddToHook{env/theorem/after}{\addvspace{3pt}}
\AddToHook{env/example/before}{\addvspace{3pt}}
\AddToHook{env/example/after}{\addvspace{3pt}}
\AddToHook{env/definition/before}{\addvspace{3pt}}
\AddToHook{env/definition/after}{\addvspace{3pt}}
\AddToHook{env/proposition/before}{\addvspace{3pt}}
\AddToHook{env/proposition/after}{\addvspace{3pt}}
\AddToHook{env/lemma/before}{\addvspace{3pt}}
\AddToHook{env/lemma/after}{\addvspace{3pt}}
\AddToHook{env/remark/before}{\addvspace{3pt}}
\AddToHook{env/remark/after}{\addvspace{3pt}}
\AddToHook{env/corollary/before}{\addvspace{3pt}}
\AddToHook{env/corollary/after}{\addvspace{3pt}}

\numberwithin{equation}{section}
\numberwithin{figure}{section}
\newlist{propenum}{enumerate}{1}
\setlist[propenum]{label=(\arabic*), ref=\theproposition~(\arabic*)}
\crefalias{propenumi}{proposition}
\newlist{lemenum}{enumerate}{1}
\setlist[lemenum]{label=(\arabic*), ref=\thelemma~(\arabic*)}
\crefalias{lemenumi}{lemma}
\newlist{thmenum}{enumerate}{1}
\setlist[thmenum]{label=(\arabic*), ref=\thetheorem~(\arabic*)}
\crefalias{thmenumi}{theorem}
\Crefformat{enumi}{#2\textup{(#1)}#3}
\newtheorem{assertion}[theorem]{Assertion}
\AddToHook{env/assertion/before}{\addvspace{3pt}}
\AddToHook{env/assertion/after}{\addvspace{3pt}}
\AddToHook{env/assertion/begin}{\crefalias{theorem}{assertion}}
\newtheorem{claim}[theorem]{Claim}
\AddToHook{env/claim/before}{\addvspace{3pt}}
\AddToHook{env/claim/after}{\addvspace{3pt}}
\AddToHook{env/claim/begin}{\crefalias{theorem}{claim}}

\usepackage[url=true,backend=biber,style=ext-alphabetic,hyperref=true,giveninits=true,maxbibnames=99]{biblatex}
\usepackage{fontspec}
\usepackage[frozencache]{minted}
\newmintinline[lean]{lean4}{fontsize=\footnotesize,bgcolor=white}
\usepackage{tcolorbox}
\tcbuselibrary{listings, minted, skins}

\tcbset{listing engine=minted}
\newtcblisting{leancode}[1][]{listing only, minted language=lean4, minted style=tango,
  colback=white!87!Emerald, enhanced, frame hidden, minted options={
  fontsize=\footnotesize, tabsize=2, breaklines, autogobble,mathescape},top=0pt,bottom=0pt,left=0pt,right=0pt,#1}

\author{Yijun Yuan\orcidlink{0000-0001-6571-6980}}
\address{Institute for Theoretical Sciences, Westlake University, No. 600 Dunyu Road, Sandun town, Xihu district, Hangzhou, Zhejiang Province, 310030, China}
\email{941201yuan@gmail.com}
\urladdr{https://yijunyuan.github.io/}

\AtBeginDocument{%
}

\begin{document}
\frontmatter
%\dedicatory{dedication}
\setttsize{\footnotesize}
\begin{abstract}
    In this article, we present a formalization of spherically complete spaces, a fundamental notion in non-Archimedean functional analysis, using the Lean theorem prover (v4.31.0), building over Mathlib. This work includes the equivalent definitions of spherically complete spaces, their basic properties, examples and non-examples such as the field $\bfC_p$ of $p$-adic complex numbers. As applications, we formalize the notion of Birkhoff–James orthogonality, the Hahn–Banach extension theorem and the spherical completion for non-Archimedean Banach spaces.
    
    URL of code: \url{https://github.com/YijunYuan/SphericalCompleteness/tree/paper}
\end{abstract}

\maketitle

\tableofcontents
\mainmatter
\section{Introduction}

\subsection{Functional analysis in the non-Archimedean setting}
Broadly speaking, functional analysis studies vector spaces endowed with compatible algebraic and topological structures, together with the continuous linear maps between them. Historically, spaces of functions provided the main examples, but the subject also treats abstract Banach spaces, Hilbert spaces and more general topological vector spaces. Over the classical fields $\bfR$ and $\bfC$, one of its fundamental results is the Hahn–Banach theorem: every continuous linear functional defined on a subspace of a normed space extends to the whole space without increasing its norm.

These foundational questions continue to hold significance when $\bfR$ or $\bfC$ is replaced by a \textbf{non-Archimedean} field, notably the field $\bfQ_p$ of $p$-adic numbers and its completed algebraic closure $\bfC_p$. Such fields play a pivotal role in modern number theory and $p$-adic geometry. In this context, the absolute value obeys the strong ultrametric inequality
$$\lVert x+y\rVert\leq \max(\lVert x\rVert, \lVert y\rVert).$$
Because this condition is strictly stronger than the standard triangle inequality, the geometric behavior of balls diverges significantly from the classical intuition, resulting in the failure or necessary modification of many standard functional analytic theorems. We refer to \Cref{sec:background} for a more detailed treatment of these foundational notions and their formalization in Lean theorem prover and its mathematical library Mathlib.

A metric space is termed \textbf{spherically complete} if every decreasing sequence of closed balls has a non-empty intersection. While this property is trivially satisfied by finite-dimensional Banach spaces over $\bfR$ or $\bfC$, it emerges as an indispensable condition in non-Archimedean functional analysis:
\begin{enumerate}
    \item The non-Archimedean Hahn–Banach extension theorem holds true for spherically complete non-Archimedean normed vector spaces (cf. \cite[Theorem 4.8]{rooijNonArchimedeanFunctionalAnalysis1978}).
    \item Not every subspace of a non-Archimedean normed vector space has an orthogonal complement. However, if the subspace is spherically complete, then it does have an orthogonal complement (cf. \cite[Corollary 4.7]{rooijNonArchimedeanFunctionalAnalysis1978}).
\end{enumerate}
Beyond functional analysis, the concept of spherical completeness is deeply intertwined with $p$-adic geometry (e.g. \cite[Proposition 1.4.4]{berkovicSpectralTheoryAnalytic2012}, \cite[Lemma 5.2.7]{scholzeModulipdivisibleGroups2013}) and transcendental number theory (e.g. \cite{kedlayaPowerSeriesPAdic2001, kedlayaAlgebraicityGeneralizedPower2017, wang2025padictranscendencesumk1inftyp1pk}).

On the other hand, spherical completeness is a geometric rather than topological condition (cf. \cite[Exercise 2.F]{rooijNonArchimedeanFunctionalAnalysis1978}). This forces operations such as taking the ``spherical completion'' (cf. \Cref{sec:27285}) of ultrametric spaces to be defined in a non-canonical way, in contrast to the canonical completion of metric spaces.

Given the central role of spherical completeness in these areas, and subtleties such as the non-canonicity of spherical completions, it is desirable to formalize the theory and its main consequences in a proof assistant, so that the relevant definitions, hypotheses and proofs can be checked with a level of precision that is difficult to achieve in informal mathematics.

\subsection{Lean, Mathlib, and related works}
\subsubsection*{Overview of Lean and Mathlib}
The Lean theorem prover (cf. \cite{demouraLeanTheoremProver2015}), or Lean for short, is an open-source interactive theorem prover and programming language based on dependent type theory. It was created by Leonardo de Moura in 2013, and the development of Lean 4 began in April 2018 \cite{leanFROAbout2026}. The current generation, Lean 4, is designed not only as a proof assistant but also as an efficient functional programming language with metaprogramming facilities, which makes it possible to write tactics, automation and domain-specific infrastructure within Lean itself. The system description of Lean 4 appeared in 2021 \cite{mouraLean4TheoremProver2021}, and Lean 4.0 was released in September 2023 \cite{leanFROAbout2026}.

The main mathematical library for Lean is Mathlib (cf. \cite{mathlib2020}). It is a community-maintained library covering large parts of algebra, topology, analysis, geometry, number theory, category theory and combinatorics, together with extensive general-purpose infrastructure for rewriting, simplification, typeclass inference and automation. As of July 2026, the official Mathlib statistics page lists 133,810 definitions, 281,979 theorems and 772 contributors \cite{mathlibStatistics2026}. These numbers make Mathlib one of the largest libraries of formalized mathematics currently available.

\subsubsection*{$p$-adic mathematics and functional analysis in Mathlib}
Currently, Mathlib already contains substantial infrastructure relevant to this project. We give the following non-exhaustive list of relevant works:
\begin{enumerate}
    \item Robert Y. Lewis formalized the field of $p$-adic numbers $\bfQ_p$ and the corresponding Hensel's lemma (cf. \cite{10.1145/3293880.3294089}).
    \item María Inés de Frutos-Fernández formalized the field $\bfC_p$ of the completed algebraic closure of $\bfQ_p$ (cf. \cite{defrutosfernandez:LIPIcs.ITP.2023.13}).
    \item Jiedong Jiang formalized Krasner's lemma (cf. \cite{jjdkrasner}) and the continuity of roots (cf. \cite{jjdcpalgclosed}). As a result, he was able to prove that $\bfC_p$ is algebraically closed (cf. \cite{jjdcpalgclosed}).
    \item Yury Kudryashov formalized the Hahn–Banach theorem for spaces over $\bfR$ (cf. \cite{KudryashovHahnBanach}).
    \item  Frédéric Dupuis, Robert Y. Lewis and Heather Macbeth formalized several fundamental results (e.g. the Fréchet–Riesz representation theorem and the spectral theorem for compact self-adjoint operators) in functional analysis for spaces over $\bfR$ or $\bfC$ (cf. \cite{dupuis_et_al:LIPIcs.ITP.2022.10}).
\end{enumerate}
%The aforementioned works in Lean 3 have all been ported to Lean 4.
\begin{remark}
    There are also related formalizations in other proof assistants: 
    \begin{enumerate}
        \item $p$-adic numbers have been formalized in Coq/Rocq by Álvaro Pelayo, Vladimir Voevodsky and Michael A. Warren (cf. \cite{pelayoUnivalentFormalizationPadic2015,UniMath}); in HOL Light by John Harrison (cf. \cite{HarrisonPadic}); and in Isabelle by Aaron Crighton (cf. \cite{Padic_Ints-AFP}).
        \item Functional analysis --- counting the Hahn–Banach theorem --- has been formalized in Coq/Rocq by Marie Kerjean and Assia Mahboubi (cf. \cite{coqhb}); in Isabelle/HOL by Gertrud Bauer (cf. \cite{isabellehb2,isabellehb1}); and in Mizar by Bogdan Nowak and Andrzej Trybulec (cf. \cite{mizarhb}).
    \end{enumerate}
\end{remark}

\subsection{Contribution and structure of this work}
In this work, we present a comprehensive formalization of results related to spherical completeness in the Lean theorem prover 4.31.0, building over Mathlib(\href{https://github.com/leanprover-community/mathlib4/tree/fabf563a7c95a166b8d7b6efca11c8b4dc9d911f}{commit fabf563})\footnote{This work was first developed in Lean 4.28.0 with Mathlib (\href{https://github.com/leanprover-community/mathlib4/tree/8f9d9cff6bd728b17a24e163c9402775d9e6a365}{commit 8f9d9cf}) and then ported to Lean 4.31.0 with Mathlib (\href{https://github.com/leanprover-community/mathlib4/tree/fabf563a7c95a166b8d7b6efca11c8b4dc9d911f}{commit fabf563}). Unless a version is given, Lean/Mathlib refers to Lean 4.31.0 and commit fabf563.}. The paper is organized as follows:
In \Cref{sec:background}, we give a brief overview of the relevant background in functional analysis, $p$-adic mathematics and its formalization in Mathlib, which provides the reader with the necessary context. In \Cref{sec:17834}, we introduce the formalization of spherical completeness, its equivalent definitions and its relationship with other topological conditions. In \Cref{sec:62478}, we formalize a criterion for non-spherical-completeness and apply it to show that $\bfC_p$ is not spherically complete. In particular, we identify and rigorously repair a gap in a classical result of Schikhof (cf. \Cref{ass:53575}). In \Cref{sec:9468}, we present our formalization of the notion of Birkhoff–James orthogonality for non-Archimedean normed spaces and establish the existence of orthogonal complements for spherically complete subspaces. In \Cref{sec:45266}, we demonstrate the formalized Hahn–Banach extension theorem for non-Archimedean normed spaces. Finally, in \Cref{sec:27285}, we present the formalization of the spherical completion of non-Archimedean normed spaces and its analytical properties.

We hope that this formalization will serve as a foundation for further developments in the field.

\subsection{Future work}
There are numerous potential directions for further formalization in non-Archimedean functional analysis and $p$-adic geometry that are related to spherically complete spaces. Among them, we highlight the following two:
\subsubsection{Spherical completeness beyond metric spaces}
Spherical completeness not only serves as an important condition in non-Archimedean functional analysis, but is also deeply connected to the model theory of non-Archimedean valued fields (cf. \cite[Section 4.2]{vandenDries2014}, \cite{williamsSphericallyCompleteModels2023}). However, the valuation does not induce a distance function when its rank is greater than 1. It might not even be metrizable. To cover these cases, one may need to generalize the notion of metric space, where the distance function takes values in a totally ordered abelian group instead of the (non-negative) real numbers. The spherical completeness condition on these generalized metric spaces is not well-studied in the literature, and it would be interesting to develop the theory and formalize it in Lean.

\subsubsection{Reflexivity}
Although the importance of spherical completeness in non-Archimedean functional analysis is well justified, this does not render non-spherically-complete spaces useless. On the contrary, spherical completeness is of no importance, or even detrimental, in Banach algebra or measure theory. For example, as was pointed out to the author by David Loeffler, there is no infinite-dimensional reflexive\footnote{A normed vector space $M$ is reflexive if the canonical embedding from $M$ to its double dual is surjective. Reflexive spaces admit many nice analytic properties. See, for example, \cite[Section 15]{schneiderNonarchimedeanFunctionalAnalysis2002a} and \cite{VANDERPUT1970341}.} Banach space over a spherically complete ultrametric normed field (cf. \cite[Theorem 4.16]{rooijNonArchimedeanFunctionalAnalysis1978}). The formalization of reflexivity and its applications in non-Archimedean functional analysis is an interesting topic and a natural continuation of this work.
\iffalse
\begin{enumerate}
    \item The spherical completeness condition can also be formulated in the context of valued fields, which coincides with the spherical completeness we define here when the valuation is of rank $1$ (cf. \cite[Definition 6.6]{barriacomicheoSummaryNonArchimedeanValued2018}). When considering valued fields, spherical completeness is related to many important topics in valuation theory, number theory and model theory, such as maximally complete valued fields (cf. \cite[Definition 6.8]{barriacomicheoSummaryNonArchimedeanValued2018}), immediate extensions (cf. \cite[Section 5]{poonenMAXIMALLYCOMPLETEFIELDS1993}) and pseudo-Cauchy sequences (cf. \cite[Definition 6.1]{barriacomicheoSummaryNonArchimedeanValued2018}) and even tilting of perfectoid fields\footnote{It can be proved that the tilt of a spherically complete perfectoid field is still spherically complete.}. It would be interesting to formalize these theories in Lean.
    \item An important class of spherically complete spaces (valued fields) is given by the fields of Hahn series and their mixed-characteristic analogues. These fields are useful for transcendental number theory over local fields and are related to the construction of big period rings (cf. \cite[page 333]{kedlayaPowerSeriesPAdic2001}) in $p$-adic Hodge theory. These constructions can be formalized in Lean.
\end{enumerate}
\fi

\subsection*{Code Availability}
The GitHub repository for this project (\url{https://github.com/YijunYuan/SphericalCompleteness}) provides two branches: \texttt{master}, which contains the latest development version, and \texttt{paper}, which preserves the exact code used in this paper for reproducibility.

Every code snippet highlighted with a colored background in this paper links directly to its corresponding implementation in the \texttt{paper} branch of the GitHub repository.

\subsection*{Declaration of AI usage}
The author used GPT 5.3 from OpenAI and Claude Opus 4.8 from Anthropic in this work in the following aspects:
\begin{enumerate}
    \item The docstrings of the Lean code are partially generated by GPT-5.3 and Claude Opus 4.8, which we fed with Mathlib's official guideline for \href{https://leanprover-community.github.io/contribute/doc.html}{documentation style} as part of the prompt.
    \item Although all the code in the initial version is completely human-written, we used Claude Opus 4.8 to simplify the code and to provide suggestions for refactoring. Mathlib's official guidelines for \href{https://leanprover-community.github.io/contribute/naming.html}{naming conventions} and \href{https://leanprover-community.github.io/contribute/naming.html}{code style} were also provided as part of the prompt.
\end{enumerate}
The author takes full responsibility for the correctness of the code and the content of this article.

\subsection*{Acknowledgements}
The author would like to thank Jiedong Jiang for helpful suggestions on the design of the formalization, and the anonymous referees for their careful reading of the manuscript and their valuable comments, which significantly improved the quality of this article and the formalized code. The research is partially supported by the National Key R\&D Program of China (Grant No. 2024YFA1014000).

\section{Preliminaries}\label{sec:background}
This section recalls background material used throughout the paper. The aim is not to give a comprehensive introduction to non-Archimedean analysis, but to fix terminology and to make clear which parts of the mathematical setting are already available in Mathlib.

\subsection{Normed vector spaces over normed fields}
\begin{definition}\leavevmode
    \begin{enumerate}
        \item A \textbf{pseudo-metric space} $(X,d)$ is a set $X$ equipped with a distance function $d\colon X\times X\to\bfR_{\geq 0}$ satisfying the following conditions:
        \begin{itemize}
            \item $d(x,x)=0$ for all $x\in X$,
            \item $d(x,y) = d(y,x)$ for all $x,y \in X$,
            \item $d(x,z) \leq d(x,y) + d(y,z)$ for all $x,y,z \in X$.
        \end{itemize}
        It is said to be a \textbf{metric space} if in addition $d(x,y) = 0$ implies $x=y$. This is formalized as a \lean{class} \href{https://github.com/leanprover-community/mathlib4/blob/fabf563a7c95a166b8d7b6efca11c8b4dc9d911f/Mathlib/Topology/MetricSpace/Pseudo/Defs.lean#L122-L152}{\lean{PseudoMetricSpace}} (resp. \href{https://github.com/leanprover-community/mathlib4/blob/fabf563a7c95a166b8d7b6efca11c8b4dc9d911f/Mathlib/Topology/MetricSpace/Defs.lean#L54-L71}{\lean{MetricSpace}}) in Mathlib.
        \item A \textbf{normed field} $K$ is a field equipped with a metric $d$ that is induced by a multiplicative norm function $\lVert\cdot\rVert\colon K\to\bfR_{\geq 0}$ via $d(x,y)=\lVert x-y\rVert$.
    \begin{leancode}[hyperurl={https://github.com/leanprover-community/mathlib4/blob/fabf563a7c95a166b8d7b6efca11c8b4dc9d911f/Mathlib/Analysis/Normed/Field/Basic.lean\#L150-L155}]
class NormedField (α : Type*) extends Norm α, Field α, MetricSpace α where
  dist_eq : ∀ x y, dist x y = norm (-x + y)
  protected norm_mul : ∀ a b, norm (a * b) = norm a * norm b
    \end{leancode}
    \noindent It is said to be \textbf{non-trivial} if there exists an element $a\in K$ such that $\lVert a\rVert>1$. This is formalized as an extended \lean{class} \href{https://github.com/leanprover-community/mathlib4/blob/fabf563a7c95a166b8d7b6efca11c8b4dc9d911f/Mathlib/Analysis/Normed/Field/Basic.lean#L157-L162}{\lean{NontriviallyNormedField}} over \lean{NormedField} with single predicate.
    \item A \textbf{semi-normed vector space} (resp. \textbf{normed vector space}) $E$ over a normed field $K$ is a vector space equipped with a pseudo-metric (resp. metric) $d$ that is induced by a semi-norm (resp. norm) function $\lVert\cdot\rVert\colon E\to\bfR_{\geq 0}$ via $d(x,y)=\lVert x-y\rVert$, satisfying that for any $c\in K$ and $v\in E$, we have $\lVert c\cdot v\rVert=\lVert c\rVert\cdot\lVert v\rVert$. This is formalized as a \lean{class} \href{https://github.com/leanprover-community/mathlib4/blob/fabf563a7c95a166b8d7b6efca11c8b4dc9d911f/Mathlib/Analysis/Normed/Module/Basic.lean#L32-L42}{\lean{NormedSpace}}, which should be used in cooperation with \href{https://github.com/leanprover-community/mathlib4/blob/fabf563a7c95a166b8d7b6efca11c8b4dc9d911f/Mathlib/Analysis/Normed/Group/Defs.lean#L223-L229}{\lean{SeminormedAddCommGroup}} (resp. \href{https://github.com/leanprover-community/mathlib4/blob/fabf563a7c95a166b8d7b6efca11c8b4dc9d911f/Mathlib/Analysis/Normed/Group/Defs.lean#L245-L250}{\lean{NormedAddCommGroup}}) in Mathlib.
    \end{enumerate}
\end{definition}
\begin{remark}
    The \lean{class} \lean{NormedSpace} in Mathlib only requires \lean{SeminormedAddCommGroup} as an implicit argument instead of \lean{NormedAddCommGroup}. As a result, \lean{NormedSpace} itself should be understood as a semi-normed space in the usual literature. In this work, part of the results around normed vector spaces are formalized in the semi-normed setting for generality, but we state most of them for normed vector spaces in this article for simplicity.
\end{remark}
\begin{remark}
    Although not used in this article, it is worth mentioning that \textbf{Banach spaces}, i.e., the complete normed vector spaces, are not represented by a standalone \lean{class}, but a combination of \lean{NormedSpace} and \href{https://github.com/leanprover-community/mathlib4/blob/fabf563a7c95a166b8d7b6efca11c8b4dc9d911f/Mathlib/Topology/UniformSpace/Cauchy.lean#L366-L370}{\lean{CompleteSpace}}.
\end{remark}

One can consider the subspaces and quotient spaces of a semi-normed vector space $E$ over a normed field $K$. Obviously, a subspace of $E$ (\href{https://github.com/leanprover-community/mathlib4/blob/fabf563a7c95a166b8d7b6efca11c8b4dc9d911f/Mathlib/Algebra/Module/Submodule/Defs.lean#L38-L42}{\lean{Submodule K E}} in Mathlib) automatically inherits the norm from $E$. On the other hand, the quotient space $E/D$ by a subspace $D\subseteq E$ carries the quotient semi-norm
$$\lVert x+D\rVert\coloneqq\inf_{y\in D}\lVert x-y\rVert.$$
This is realized in Mathlib as an \lean{instance} \href{https://github.com/leanprover-community/mathlib4/blob/fabf563a7c95a166b8d7b6efca11c8b4dc9d911f/Mathlib/Analysis/Normed/Group/Quotient.lean#L462-L464}{\lean{Submodule.Quotient.normedSpace}}, and the formula is reflected in the lemma \href{https://github.com/leanprover-community/mathlib4/blob/fabf563a7c95a166b8d7b6efca11c8b4dc9d911f/Mathlib/Analysis/Normed/Group/Quotient.lean#L160-L168}{\lean{QuotientGroup.norm_mk}}.

For (semi-)normed vector spaces $E,F$ over $K$, any continuous linear map $T\colon E\lto F$ can be endowed with the \textbf{operator norm}, which is defined as the infimum of the non-negative real numbers $c$ such that $\lVert T(x)\rVert\leq c\cdot \lVert x\rVert$ for all $x\in E$. This is formalized in Mathlib under the name \href{https://github.com/leanprover-community/mathlib4/blob/fabf563a7c95a166b8d7b6efca11c8b4dc9d911f/Mathlib/Analysis/Normed/Operator/Basic.lean#L172-L174}{\lean{ContinuousLinearMap.opNorm}}, and one can use the notation \lean{‖T‖} to denote the operator norm of $T$. The set of all continuous linear maps from $E$ to $F$, denoted by $\calL(E,F)$, forms a (semi-)normed vector spaces over $K$ under the operator norm. This is formalized in Mathlib under the name \href{https://github.com/leanprover-community/mathlib4/blob/fabf563a7c95a166b8d7b6efca11c8b4dc9d911f/Mathlib/Analysis/Normed/Operator/Basic.lean#L390-L392}{\lean{ContinuousLinearMap.toNormedSpace}}.

%Orthogonality also needs care outside Hilbert spaces. In Hilbert spaces over $\bfR$ or $\bfC$, orthogonality is defined by an inner product. General normed spaces have no such structure. The replacement used later is metric: an element $x$ is orthogonal to a subspace $D$ if $\lVert x\rVert$ is exactly its distance from $D$, or equivalently if adding an element of $D$ cannot decrease its norm. This is a special case of Birkhoff–James orthogonality.

\subsection{Non-archimedean property and number-theoretic preliminaries}
\begin{definition}
    A metric space $(X,d)$ is \textbf{non-Archimedean} or \textbf{ultrametric} if for any $x,y,z\in X$, we have
    $$d(x,z)\leq \max\left(d(x,y),d(y,z)\right).$$
    This is defined as a \lean{class} \href{https://github.com/leanprover-community/mathlib4/blob/fabf563a7c95a166b8d7b6efca11c8b4dc9d911f/Mathlib/Topology/MetricSpace/Ultra/Basic.lean#L44-L47}{\lean{IsUltrametricDist}} in Mathlib.
\end{definition}

Ultrametric spaces carry many interesting properties that are not shared by general metric spaces. For example, a frequently used result in this work is that any two closed balls in an ultrametric space are either disjoint or one is contained in the other. This is formalized in Mathlib under the name \href{https://github.com/leanprover-community/mathlib4/blob/fabf563a7c95a166b8d7b6efca11c8b4dc9d911f/Mathlib/Topology/MetricSpace/Ultra/Basic.lean#L103-L113}{\lean{IsUltrametricDist.closedBall_subset_trichotomy}}.
\begin{leancode}[hyperurl={https://github.com/leanprover-community/mathlib4/blob/fabf563a7c95a166b8d7b6efca11c8b4dc9d911f/Mathlib/Topology/MetricSpace/Ultra/Basic.lean\#L104-L105}]
closedBall x r ⊆ closedBall y s ∨ closedBall y s ⊆ closedBall x r ∨ 
Disjoint (closedBall x r) (closedBall y s)
\end{leancode}

One of the most important motivating examples of ultrametric spaces are the $p$-adic fields.
\begin{definition}
    Let $p$ be a prime number.
    \begin{enumerate}
        \item The \textbf{$p$-adic valuation} of an integer $n\in\bfZ$ is defined as
        $$v_p(n)\coloneqq \max\{k\in\bfN \colon p^k \mid n\}\in \bfN\cup\{+\infty\},$$
        where we set $v_p(0)\coloneqq +\infty$.
        \item For a rational number $q=a/b\in\bfQ$ with $a,b\in\bfZ$, the \textbf{$p$-adic norm} of $q$ is defined as
        $$\lVert q\rVert_p\coloneqq p^{-v_p(a)+v_p(b)}.$$
        When $q=0$, we regard $\lVert 0\rVert_p$ as $0$.
        \item The field $\bfQ_p$ of $p$-adic rational numbers is defined to be the completion of $\bfQ$ with respect to the $p$-adic norm. The field $\bfQ_p$ is denoted by \href{https://github.com/leanprover-community/mathlib4/blob/fabf563a7c95a166b8d7b6efca11c8b4dc9d911f/Mathlib/NumberTheory/Padics/PadicNumbers.lean#L527-L534}{\lean{ℚ_[p]}} in Mathlib.
    \end{enumerate}
\end{definition}
It is easy to verify, and is already formalized in Mathlib as \lean{instance} \href{https://github.com/leanprover-community/mathlib4/blob/fabf563a7c95a166b8d7b6efca11c8b4dc9d911f/Mathlib/NumberTheory/Padics/PadicNumbers.lean#L779-L780}{\lean{Padic.instIsUltrametricDist}}, that the $p$-adic norm on $\bfQ_p$ is ultrametric.

The $p$-adic norm extends uniquely to the algebraic closure $\bfQ^{\opn{alg}}_p$ of $\bfQ_p$ (\href{https://github.com/leanprover-community/mathlib4/blob/fabf563a7c95a166b8d7b6efca11c8b4dc9d911f/Mathlib/NumberTheory/Padics/Complex.lean#L53-L54}{\lean{PadicAlgCl p}} in Mathlib). However, a non-trivial fact is that $\bfQ^{\opn{alg}}_p$ is not complete with respect to the extended $p$-adic norm (cf. \cite[Section 1.4]{robertCoursePadicAnalysis2000}). This allows us to enlarge $\bfQ^{\opn{alg}}_p$ to its completion, which is denoted by $\bfC_p$ (\href{https://github.com/leanprover-community/mathlib4/blob/fabf563a7c95a166b8d7b6efca11c8b4dc9d911f/Mathlib/NumberTheory/Padics/Complex.lean#L135-L140}{\lean{ℂ_[p]}} in Mathlib). Basic properties of $\bfC_p$ are already formalized in Mathlib, including the fact that the extended $p$-adic norm on $\bfC_p$ is ultrametric.

\begin{remark}
    At the moment when the initial version of this project was finished with Lean 4.28 over Mathlib (\href{https://github.com/leanprover-community/mathlib4/tree/8f9d9cff6bd728b17a24e163c9402775d9e6a365}{Commit 8f9d9cf}), there was a design flaw in the definition of the norm on $\bfC_p$ in Mathlib. The completion process automatically endows $\bfC_p$ with the extended norm, but another (mathematically identical, yet not Lean-semantically identical) norm, coming from the $p$-adic valuation, was written manually. This caused a diamond issue in the typeclass inference of Lean. It was fixed in a later version of Mathlib by Jiedong Jiang (cf. \cite{jjddiamond}).
\end{remark}

\subsection{Conventions and notations}
In the rest of this article, we adopt the following conventions and notations:
\begin{itemize}
    \item $p$ will always denote a prime number.
    \item $(\bbK,\lVert\cdot \rVert)$ will be a non-trivially ultrametric normed field. Unless otherwise specified, any other fields appearing in this article are not assumed to be ultrametric.
    \item By a \textbf{normed space over $\bbK$}, we mean an ultrametric normed vector space over $\bbK$.
    \item By a \textbf{ball} in a metric space, we always mean a closed ball with non-negative radius.
    %\item For any normed spaces $E,F$ over $\bbK$, we denote by $\calL(E,F)$ the space of continuous linear maps from $E$ to $F$, equipped with the operator norm.
    \item For a vector space $E$ over $\bbK$ and an element $a\in E$, we denote by $[a]$ the subspace of $E$ generated by $a$.
    \item For a metric space $(X,d)$, an element $x\in X$ and a subset $A\subseteq X$, we denote by $d(x,A)$ the distance from $x$ to $A$, i.e. $d(x,A)\coloneqq \inf_{y\in A} d(x,y)$. In Mathlib, this is named as \href{https://github.com/leanprover-community/mathlib4/blob/fabf563a7c95a166b8d7b6efca11c8b4dc9d911f/Mathlib/Topology/MetricSpace/HausdorffDistance.lean#L572-L574}{\lean{Metric.infDist}}.
\end{itemize}
  
\begin{remark}
    Some results in this article for metric spaces (resp. normed spaces over $\bbK$) can be slightly generalized for pseudo-metric spaces (resp. semi-normed spaces over $\bbK$), as we do in the formalization. However, we only state and prove the results for metric spaces (resp. normed spaces over $\bbK$) in this article for simplicity.
\end{remark}

\section{Basic properties of spherically complete spaces}\label{sec:17834}
\subsection{Definitions and basic operations}\label{sec:18000}
As mentioned in the introduction, a metric space is said to be \textbf{spherically complete} if every decreasing sequence of closed balls has a non-empty intersection. This is formalized as follows:
\begin{leancode}[hyperurl={https://github.com/YijunYuan/SphericalCompleteness/blob/717759d2b22a5720e74c14160a8385ac17f32abc/SphericalCompleteness/Defs.lean\#L25-L39}]
class SphericallyCompleteSpace (α : Type*) [PseudoMetricSpace α] : Prop where
  nonempty_iInter_of_antitone : ∀ ⦃ci : ℕ → α⦄, ∀ ⦃ri : ℕ → ℝ≥0⦄,
    Antitone (fun i => closedBall (ci i) (ri i)) →
      (⋂ i, closedBall (ci i) (ri i)).Nonempty
\end{leancode}
Here, \lean{Antitone} means that the sequence of closed balls is decreasing.
\begin{remark}
    Note that \href{https://github.com/leanprover-community/mathlib4/blob/fabf563a7c95a166b8d7b6efca11c8b4dc9d911f/Mathlib/Topology/MetricSpace/Pseudo/Defs.lean#L417-L419}{\lean{Metric.closedBall}} takes an arbitrary real number as its radius. To ensure that the balls are non-empty, i.e. that the radius is non-negative, we use \href{https://github.com/leanprover-community/mathlib4/blob/fabf563a7c95a166b8d7b6efca11c8b4dc9d911f/Mathlib/Data/NNReal/Defs.lean#L57-L62}{\lean{ℝ≥0}}, the type of non-negative real numbers in Mathlib, instead of \href{https://github.com/leanprover-community/mathlib4/blob/fabf563a7c95a166b8d7b6efca11c8b4dc9d911f/Mathlib/Data/Real/Basic.lean#L33-L41}{\lean{Real}}. As a result, an implicit coercion from \lean{ℝ≥0} to \lean{Real} is applied when calling \lean{closedBall}.
\end{remark}

\begin{example}
    The most basic examples of spherically complete spaces are the fields $\bfR$ or $\bfC$ with the standard norm and the field $\bfQ_p$ of $p$-adic numbers with the $p$-adic norm. We will show in \Cref{sec:18500} that the spherical completeness of these fields can be proved in a unified way.
\end{example}

In practice, it is often more convenient to work with sequences of balls with decreasing radius, which is equivalent to the above definition:
\begin{theorem}[{cf. \cite[Lemma 2.3]{rooijNonArchimedeanFunctionalAnalysis1978}}]\label{thm:5601}
    For an ultrametric space $(X,d)$, the following are equivalent:
    \begin{enumerate}
        \item $X$ is spherically complete;
              %\item Every sequence of balls with decreasing radius has a non-empty intersection.
        \item Every sequence of balls with strictly decreasing radius has a non-empty intersection.
        \item Any collection of balls with pairwise non-empty intersection has a non-empty intersection.
    \end{enumerate}
\end{theorem}
This theorem is formalized in Lean as follows:
\begin{leancode}[hyperurl={https://github.com/YijunYuan/SphericalCompleteness/blob/717759d2b22a5720e74c14160a8385ac17f32abc/SphericalCompleteness/Basic.lean\#L99-L141}]
theorem SphericallyCompleteSpace.iff_strictAnti_radius : -- (1) $\Leftrightarrow$ (2)
    SphericallyCompleteSpace α ↔
    ∀ ⦃ci : ℕ → α⦄, ∀ ⦃ri : ℕ → NNReal⦄, StrictAnti ri →
      Antitone (fun i ↦ closedBall (ci i) (ri i)) → 
      (⋂ i, closedBall (ci i) (ri i)).Nonempty
\end{leancode}
\begin{leancode}[hyperurl={https://github.com/YijunYuan/SphericalCompleteness/blob/717759d2b22a5720e74c14160a8385ac17f32abc/SphericalCompleteness/Basic.lean\#L262-L326}]
theorem SphericallyCompleteSpace.iff_pairwise_inter_nonempty : -- (1) $\Leftrightarrow$ (3)
    SphericallyCompleteSpace α ↔ 
    (∀ S : Set (α × NNReal), S.Nonempty →
       (∀ w1 w2 : ↑S, (closedBall w1.val.1 w1.val.2 ∩ closedBall w2.val.1 w2.val.2).Nonempty) →
       (⋂ w : ↑S, closedBall w.val.1 w.val.2).Nonempty)
\end{leancode}
\iffalse
\begin{remark}
    Lean does not accept ``anonymous'' ball, i.e. ball without explicit center and radius. Therefore, to express a collection of balls, we use a set of pairs of center and radius.
\end{remark}

\fi

We conclude this section by recalling the following properties of spherically complete spaces, which are easy to prove and formalize, and useful in practice:
\begin{proposition}\label{prop:1411}\leavevmode
    \begin{enumerate}
        \item The product of finitely many spherically complete metric spaces is again spherically complete.
        \item Spherical completeness is preserved under isometric isomorphisms.
        \item Any quotient of a spherically complete normed space over $\bbK$ is spherically complete.
    \end{enumerate}
\end{proposition}
This theorem is formalized in this project as follows:
\begin{leancode}[hyperurl={https://github.com/YijunYuan/SphericalCompleteness/blob/717759d2b22a5720e74c14160a8385ac17f32abc/SphericalCompleteness/Basic.lean\#L363-L443}]
-- Product of 2 spherically complete spaces is still spherically complete
instance Prod.sphericallyCompleteSpace {E F : Type*}
    [PseudoMetricSpace E] [PseudoMetricSpace F]
    [SphericallyCompleteSpace E] [SphericallyCompleteSpace F] :
    SphericallyCompleteSpace (E × F)

-- Product of a finite family of spherically complete spaces is still spherically complete
instance Pi.sphericallyCompleteSpace {ι : Type*} [Fintype ι] {E : ι → Type*}
    [∀ i, PseudoMetricSpace (E i)] [∀ i, SphericallyCompleteSpace (E i)] :
    SphericallyCompleteSpace (∀ i, E i)
\end{leancode}
\begin{leancode}[hyperurl={https://github.com/YijunYuan/SphericalCompleteness/blob/717759d2b22a5720e74c14160a8385ac17f32abc/SphericalCompleteness/Defs.lean\#L86-L114}]
-- Spherical completeness is preserved under isometric isomorphisms
theorem SphericallyCompleteSpace.of_isometryEquiv {E F : Type*}
    [PseudoMetricSpace E] [PseudoMetricSpace F]
    [SphericallyCompleteSpace E]
    (f : E ≃ᵢ F) : SphericallyCompleteSpace F
\end{leancode}
\begin{leancode}[hyperurl={https://github.com/YijunYuan/SphericalCompleteness/blob/717759d2b22a5720e74c14160a8385ac17f32abc/SphericalCompleteness/NormedVectorSpace/Quotient.lean\#L135-L183}]
-- Any quotient of spherically complete normed spaces over 𝕂 is spherically complete
theorem Quotient.sphericallyCompleteSpace (𝕜 : Type*) [NontriviallyNormedField 𝕜]
    {E : Type*} [SeminormedAddCommGroup E] [NormedSpace 𝕜 E]
    [IsUltrametricDist E] [SphericallyCompleteSpace E]
    {F : Submodule 𝕜 E} :
    SphericallyCompleteSpace (E ⧸ F)
\end{leancode}

\subsection{Relation with other topological conditions}\label{sec:15313}
\subsubsection{Completeness}
It is a classical result that a metric space is complete if and only if every decreasing sequence of closed balls with radius tending to zero has a non-empty intersection, which is formalized as an external lemma in our project:
\begin{leancode}[hyperurl={https://github.com/YijunYuan/SphericalCompleteness/blob/717759d2b22a5720e74c14160a8385ac17f32abc/SphericalCompleteness/External/Complete.lean\#L26-L80}]
theorem completeSpace_iff_nonempty_iInter_closedBall_of_tendsto_zero
    (α : Type*) [PseudoMetricSpace α] :
    CompleteSpace α ↔
    ∀ ⦃ci : ℕ → α⦄ ⦃ri : ℕ → NNReal⦄,
      Antitone (fun i ↦ closedBall (ci i) (ri i)) →
      Filter.Tendsto ri atTop (nhds 0) →
      (⋂ i, closedBall (ci i) (ri i)).Nonempty
\end{leancode}

As an immediate corollary, every spherically complete space is complete, and this is formalized as an \lean{instance} in our project.

\begin{remark}
    As suggested by the anonymous referee, we assign this \lean{instance} a lower priority than the default one, so that Lean is less likely to prove the completeness of a space via spherical completeness, which is rarely needed in practice.
\end{remark}

\subsubsection{Properness}\label{sec:18500}
The condition of spherical completeness is fulfilled in various contexts that arise in different settings:
\begin{enumerate}
    \item A subspace of a locally compact normed vector space over a complete non-trivially normed field is spherically complete.
    \item Finite-dimensional subspaces of normed vector spaces over locally compact non-trivially normed fields are spherically complete. Note that this implies that every finite-dimensional normed space over $\bfR$, $\bfC$ and $\bfQ_p$ is spherically complete.
    \item Compact metric spaces are spherically complete.
    \item The complex upper half-plane, equipped with the hyperbolic metric, is spherically complete.
\end{enumerate}
It is possible, but very inefficient, to formalize these results individually in Lean. Instead, we observe that all these metric spaces are \textbf{proper}, i.e. every closed ball is compact. As a result, the spherical completeness of these spaces follows from Cantor's intersection theorem (cf. \cite[Theorem 3.25]{apostolMathematicalAnalysisModern1974}), i.e. the fact that the intersection of a decreasing sequence of non-empty compact sets is non-empty:
\begin{leancode}[hyperurl={https://github.com/YijunYuan/SphericalCompleteness/blob/717759d2b22a5720e74c14160a8385ac17f32abc/SphericalCompleteness/Defs.lean\#L61-L84}]
instance instSphericallyCompleteSpaceOfProperSpace (α : Type*)
    [PseudoMetricSpace α] [ProperSpace α] : SphericallyCompleteSpace α where
  nonempty_iInter_of_antitone := by
    intro ci ri hanti
    exact IsCompact.nonempty_iInter_of_sequence_nonempty_isCompact_isClosed
      (fun i ↦ closedBall (ci i) ↑(ri i)) -- chain of balls
      (fun _ ↦ hanti (by linarith)) -- the chain is descending 
      (fun h ↦ nonempty_closedBall.mpr (ri h).prop) -- balls are non-empty
      (isCompact_closedBall (ci 0) ↑(ri 0)) -- the first ball is compact
      fun _ ↦ isClosed_closedBall -- all balls are closed
\end{leancode}

\begin{remark}
    The properness of the aforementioned metric spaces is already formalized in Mathlib. With this \lean{instance}, Lean can automatically infer the spherical completeness of these spaces without extra effort.
\end{remark}

\section{Non-examples: separable spherically dense ultrametric spaces}\label{sec:62478}
Among complete valued fields, $\bfC_p$, the completion of the algebraic closure of $\bfQ_p$, is of particular importance in $p$-adic geometry. It is therefore desirable to formalize the non-spherical completeness of $\bfC_p$ in Lean.

\subsection{A second glance at a result of Schikhof}
\begin{definition}\leavevmode
    \begin{enumerate}
        \item A metric space $(X,d)$ is said to have \textbf{dense metric}, if for any ball $B$ in $X$, the set $d(B\times B)$ is dense in $[0,\opn{diam}(B)]$.
        \begin{leancode}[hyperurl={https://github.com/YijunYuan/SphericalCompleteness/blob/717759d2b22a5720e74c14160a8385ac17f32abc/SphericalCompleteness/Dense.lean\#L476-L493}]
class IsDenseMetric (α : Type*) [MetricSpace α] : Prop where
  dense_metric : ∀ z : α, ∀ r : ℝ≥0,
    closure (Set.image2 dist (closedBall z r) (closedBall z r)) 
      = Set.Icc 0 (diam (closedBall z r))
        \end{leancode}
        \item A topological space is \textbf{separable} if it has a countable dense subset.
    \end{enumerate}
\end{definition}
Schikhof provides a proof of the following assertion:
\begin{assertion}[{cf. \cite[Theorem 20.5]{schikhofUltrametricCalculusIntroduction1985}}]\label{ass:53575}
    Separable ultrametric spaces with dense metric are not spherically complete.
\end{assertion}
When formalizing this result in Lean, we find that the denseness of the metric is not enough to ensure non-spherical completeness. In particular, a metric space consisting of a single point trivially has dense metric and is separable, but it is also spherically complete, as there is only one ball in this space.
\begin{remark}
    To present this counterexample accurately, we establish \lean{instance}s for \lean{PUnit}, which is the \lean{Type} representing a single point.
    \begin{leancode}[hyperurl={https://github.com/YijunYuan/SphericalCompleteness/blob/717759d2b22a5720e74c14160a8385ac17f32abc/SphericalCompleteness/Dense.lean\#L645-L653}]
-- Single point space is a separable ultrametric space with dense metric.
example : SeparableSpace PUnit := inferInstance
example : IsDenseMetric PUnit := inferInstance
example : IsUltrametricDist PUnit := inferInstance
-- Single point space is spherically complete, which disproofs Schikhof's result.
example : SphericallyCompleteSpace PUnit := inferInstance    
    \end{leancode}
\end{remark}
\begin{remark}
    The single point space is not the only counterexample to Schikhof's result. For instance, see \cite{5132970} for a more non-trivial counterexample constructed by Torsten Schoeneberg.
\end{remark}

To isolate the issue, we revisit Schikhof's proof. Suppose, for contradiction, that $(X,d)$ is a spherically complete separable ultrametric space with dense metric. Let $\{a_0,a_1,\cdots\}$ be a countable\footnote{Note that this set can be finite. Since the proof is similar, we assume that this set has cardinality $\aleph_0$.} dense subset of $X$. Take $r_0>0$ to be the distance between two distinct points of $X$, and take a sequence of real numbers $r_1,r_2,\cdots$ such that $r_0>r_1>\cdots>r_n>\cdots>\frac{1}{2}r_0$.

Schikhof's proof depends on the following claim:
\begin{claim}\label{claim:34467}
    Let $x\in X$ and let $B\subseteq X$ be a ball with $\opn{diam}(B)>0$ or $X$ itself. For any $r'\in(0,\opn{diam}(B))$, there exists a ball $B'\subseteq B$ with $\opn{diam}(B')=r'$ such that $x\notin B'$.
\end{claim}
With this claim, one can recursively construct a decreasing sequence of balls $B_1\supseteq B_2\supseteq \cdots\supseteq B_n\supseteq\cdots$ such that $\opn{diam}(B_n)=r_n$ and $a_{n-1}\notin B_n$ for each $n\geq 1$. The spherical completeness of $X$ implies that $\bigcap_{n=1}^\infty B_n$ is non-empty and consequently contains a non-empty ball of radius $\frac{1}{2}r_0>0$. This ball is disjoint from the dense set $\{a_0,a_1,\cdots\}$, which leads to a contradiction.

The subtlety appears in the proof of \Cref{claim:34467}. Schikhof partitions $B$ into at least two disjoint balls of radius $r'$ via the relation $x\sim y \Longleftrightarrow d(x,y)\leq r'$, and one of them, say $B'$, does not contain $x$. Although this is sound, it is not clear why $\opn{diam}(B')=r'$. The denseness of the metric only ensures that there are subballs of $B$ with diameter arbitrarily close to $r'$; it does not guarantee that the distances approximating $r'$ from below can always be realized between points in the same subball, let alone that we can find such a subball that does not contain $x$.

To fix this issue, we introduce a stronger condition than the denseness of the metric:
\begin{definition}
    A metric space is \textbf{spherically dense} if for any ball $B$, its diameter equals the radius.
\end{definition}
\begin{leancode}[hyperurl={https://github.com/YijunYuan/SphericalCompleteness/blob/717759d2b22a5720e74c14160a8385ac17f32abc/SphericalCompleteness/Dense.lean\#L55-L68}]
class IsSphericallyDense (α : Type*) [PseudoMetricSpace α] : Prop where
  spherically_dense : ∀ (c : α) (r : ℝ≥0) , diam (closedBall c r) = r    
\end{leancode}
\begin{theorem}\label{thm:3701}
    Let $(X,d)$ be a separable spherically dense ultrametric space. Then \Cref{claim:34467} holds. In particular, $X$ is not spherically complete.
\end{theorem}
\begin{proof}
   Since $r'<\opn{diam}(B)$, one can take $x_1,x_2\in B$ such that $d(x_1,x_2)>r'$. Consider the balls $B(x_1,r')$ and $B(x_2,r')$. The ultrametric inequality ensures that these two balls are disjoint, and consequently at least one of them does not contain $x$. The spherical denseness of $X$ ensures that the diameter of this ball equals $r'$.
\end{proof}
In our project, \Cref{claim:34467} is formalized as follows:
\begin{leancode}[hyperurl={https://github.com/YijunYuan/SphericalCompleteness/blob/717759d2b22a5720e74c14160a8385ac17f32abc/SphericalCompleteness/Dense.lean\#L204-L240}]
lemma exists_sub_closedBall_not_mem {α : Type*}
    [PseudoMetricSpace α] [IsUltrametricDist α] [IsSphericallyDense α]
    (c₀ : α) (r₀ : ℝ≥0) (r₁ : ℝ≥0) (hr : r₁ < r₀) (z : α) :
    ∃ c₁ : α,
    closedBall c₁ r₁ ⊆ closedBall c₀ r₀ ∧
    z ∉ closedBall c₁ r₁osedBall c₁ r₁ ⊆ closedBall c₀ r₀ ∧ z ∉ closedBall c₁ r₁
\end{leancode}
And the non-spherical completeness of separable spherically dense ultrametric spaces is formalized as follows:
\begin{leancode}[hyperurl={https://github.com/YijunYuan/SphericalCompleteness/blob/717759d2b22a5720e74c14160a8385ac17f32abc/SphericalCompleteness/Dense.lean\#L358-L460}]
theorem not_sphericallyCompleteSpace_of_isSphericallyDense_separable_ultrametric
    (α : Type*) [MetricSpace α]
    [IsUltrametricDist α] [IsSphericallyDense α] [Nonempty α] [SeparableSpace α] :
    ¬ SphericallyCompleteSpace α
\end{leancode}
The following proposition shows that spherical denseness is a reasonable replacement for the dense metric condition and is satisfied by the important example $\bfC_p$:
\begin{proposition}\label{prop:7355}
    Let $(X,d)$ be an ultrametric space.
    \begin{enumerate}
        \item\label{it:58061} If $X$ is spherically dense, then it has dense metric. The converse is not true.
\begin{leancode}[hyperurl={https://github.com/YijunYuan/SphericalCompleteness/blob/717759d2b22a5720e74c14160a8385ac17f32abc/SphericalCompleteness/Dense.lean\#L496-L546}]
theorem IsDenseMetric.of_isSphericallyDense (α : Type*)
    [MetricSpace α] [IsUltrametricDist α] [IsSphericallyDense α] :
    IsDenseMetric α
\end{leancode}
        \item If $(F,\lVert\cdot\rVert)$ is a normed field, then the following are equivalent:
        \begin{enumerate}[ref=\alph*]
            \item\label{it:16856} $F$ is spherically dense;
            \item\label{it:14189} $F$ is a \textbf{densely normed field}, i.e. for any positive real numbers $r_1<r_2$, there exists $z\in F$ such that $r_1<\lVert z\rVert<r_2$.
            \item\label{it:2589} $F$ has dense metric.
        \end{enumerate}
\begin{leancode}[hyperurl={https://github.com/YijunYuan/SphericalCompleteness/blob/717759d2b22a5720e74c14160a8385ac17f32abc/SphericalCompleteness/Dense.lean\#L72-L138}]
instance instIsSphericallyDenseOfDenselyNormedField (α : Type*)
    [DenselyNormedField α] [IsUltrametricDist α] :
    IsSphericallyDense

instance instDenselyNormedFieldOfIsSphericallyDenseOfUltrametricNormedField (α : Type*)
    [NormedField α] [IsUltrametricDist α] [IsSphericallyDense α] :
    DenselyNormedField α
\end{leancode}
\begin{leancode}[hyperurl={https://github.com/YijunYuan/SphericalCompleteness/blob/717759d2b22a5720e74c14160a8385ac17f32abc/SphericalCompleteness/Dense.lean\#L548-L630}]
theorem IsDenseMetric.iff_isSphericallyDense_of_normedField (α : Type*)
    [NormedField α] [IsUltrametricDist α] :
    IsDenseMetric α ↔ IsSphericallyDense α 
\end{leancode}
    \end{enumerate}
\end{proposition}
\begin{proof}\leavevmode
    \begin{enumerate}
        \item Suppose that $(X,d)$ is spherically dense. Let $B$ be a ball in $X$ with $\opn{diam}(B)>0$, i.e. it contains at least two different points $x,y\in B$. For any real numbers $r_1<r_2$ in the interval $[0,\opn{diam}(B)]$, consider the ball $B'\coloneqq B\left(x,\frac{r_1+r_2}{2}\right)$. The ultrametric inequality ensures that $B'$ is contained in $B$. The spherical denseness of $X$ implies that $\opn{diam}(B')=\frac{r_1+r_2}{2}$. As a result, one can take $z_1,z_2\in B'\subseteq B$ such that $d(z_1,z_2)=\frac{r_1+r_2}{2}\in (r_1,r_2)$. This shows that $d(B\times B)$ is dense in $[0,\opn{diam}(B)]$.
        
        The space consisting of a single point is a counterexample to the converse.
        \item By \Cref{it:58061}, we only need to show the equivalence between \Cref{it:16856} and \Cref{it:14189}, and that \Cref{it:2589} implies \Cref{it:14189}. 
        
        Suppose that $F$ is spherically dense. For any positive real numbers $r_1<r_2$, consider the ball $B\coloneqq B\left(0,\frac{r_1+r_2}{2}\right)$. The spherical denseness of $F$ implies that $\opn{diam}(B)=\frac{r_1+r_2}{2}$. As a result, one can take $x,y\in B$ such that $d(x,y)=\lVert x-y\rVert\in \left(r_1,\frac{r_1+r_2}{2}\right)$. This proves that \Cref{it:16856} implies \Cref{it:14189}.
        
        Conversely, suppose that $F$ is a densely normed field. Let $B\coloneqq B(z,r)$ be a ball in $F$ with $r>0$. For any $\epsilon>0$, the denseness of the norm ensures that there exists $w\in F$ such that $$r-\epsilon<\lVert w\rVert=d(z,z+w)<r.$$
        This shows that $z+w\in B$ and $\opn{diam}(B)\geq r-\epsilon$. 
        The implication \Cref{it:14189}$\Rightarrow$\Cref{it:16856} follows by letting $\epsilon$ tend to zero.

        Now we suppose $F$ has dense metric. We first claim that $F$ is non-trivially valued, i.e. there exists $w\in F$ such that $\norm{w}>1$. Consider the unit ball $B_0=B(0,1)$. Since $\opn{diam}(B_0)\geq \norm{1-0}=1$, the ultrametric inequality implies that $\opn{diam}(B_0)=1$. Thus the dense metric ensures the existence of $z_1,z_2\in B_0$ such that $d(w_1,w_2)=\norm{w_1-w_2}\in(\frac{1}{2},1)$. Then taking $w=(w_1-w_2)^{-1}$ gives $\norm{w}>1$. To show that $F$ is a densely normed field, let $r_1<r_2$ be two positive real numbers. By the claim, one takes an integer $n>0$ large enough such that $R\coloneqq \norm{w}^n>r_2$. By similar reasoning as above, one knows that $\opn{diam}(B(0,R))=R$. Since $(r_1,r_2)\subseteq [0,R]$, the dense metric ensures that there exists $z_1,z_2\in B(0,R)$ such that $d(z_1,z_2)=\norm{z_1-z_2}\in(r_1,r_2)$; then $z=z_1-z_2\in F$ satisfies $r_1<\norm{z}<r_2$, which finishes the proof.
    \end{enumerate}
\end{proof}
\begin{remark}
    We do not formalize the equivalence of spherical denseness and dense metric condition for normed fields as \lean{instances}. This is because the \lean{class} \lean{IsDenseMetric} is not supposed to be used anywhere else in the project, and \lean{instance}s for it may slow down the type inference of Lean.
\end{remark}

\subsection{The non-spherical completeness of \texorpdfstring{$\bfC_p$}{Cp}}
Fix a prime number $p$:
\begin{leancode}[hyperurl={https://github.com/YijunYuan/SphericalCompleteness/blob/717759d2b22a5720e74c14160a8385ac17f32abc/SphericalCompleteness/External/PadicAlgCl.lean\#L25}]
variable (p : ℕ) [hp : Fact (Nat.Prime p)]
\end{leancode}
To apply \Cref{thm:3701} to $\bfC_p$ in Lean, one needs to provide the \lean{instance}s of \lean{DenselyNormedField} and \lean{TopologicalSpace.SeparableSpace} for $\bfC_p$, which are not yet formalized in Mathlib. In this section, we demonstrate our implementation of these two \lean{instance}s.

\subsubsection{Norm-denseness of $\bfC_p$}
Mathematically speaking, it is fairly easy to show that $\bfC_p$ is a densely normed field: for any positive real numbers $r_1<r_2$, one can take a rational number $r$ such that $r_1<\lVert p\rVert^r<r_2$. Then $p^r$ is an element of $\bfC_p$ with norm in $(r_1,r_2)$. To formalize this, we divide the proof into the following steps:
\begin{enumerate}
    \item For real numbers $b>a\geq 0$, there exists a rational number $r=\frac{u}{v}$ with $\gcd(u,v)=1$ such that $a<\lVert p\rVert^r<b$, where $p$ is viewed as an element in $\bfQ^{\opn{alg}}_p$.
    \begin{leancode}[hyperurl={https://github.com/YijunYuan/SphericalCompleteness/blob/717759d2b22a5720e74c14160a8385ac17f32abc/SphericalCompleteness/External/PadicAlgCl.lean\#L27-L56}]
lemma exists_rat_pow_p_norm_between {a b : ℝ} (ha : 0 ≤ a) (hab : a < b) : 
    ∃ c : ℚ,
    a < ‖(p : (PadicAlgCl p))‖ ^ (c : ℝ) ∧ ‖(p : (PadicAlgCl p))‖ ^ (c : ℝ) < b
    \end{leancode}
    There is nothing new in the proof: take the logarithm with base $\lVert p\rVert$ and use the density of $\bfQ$ in $\bfR$.
    \item Let $f(X)=X^v-p^u\in\bfQ^{\opn{alg}}_p[X]$, and let $z\in\bfQ^{\opn{alg}}_p$ be a root of $f$. Then one has $z^v=p^u$, and consequently $\lVert z\rVert^v=\lVert p\rVert^u$, i.e. $\lVert z\rVert=\lVert p\rVert^r$. This shows that $\bfQ^{\opn{alg}}_p$ is a densely normed field.
    \begin{leancode}[hyperurl={https://github.com/YijunYuan/SphericalCompleteness/blob/717759d2b22a5720e74c14160a8385ac17f32abc/SphericalCompleteness/External/PadicAlgCl.lean\#L58-L88}]
instance instDenselyNormedFieldPadicAlgCl : 
    DenselyNormedField (PadicAlgCl p) where
  lt_norm_lt a b ha hab := by
    rcases exists_rat_pow_p_norm_between p a b ha hab with ⟨r, hr1, hr2⟩
    let f : Polynomial (PadicAlgCl p) :=
      X ^ r.den - C ((p : PadicAlgCl p) ^ r.num) -- r.den = v, r.num = u
    ...
    rcases IsAlgClosed.exists_root f hf with ⟨z, hz⟩ -- hf : f.degree ≠ 0
    ...
    have hz' : ‖z‖ = ‖(p : PadicAlgCl p)‖ ^ (↑r : ℝ) := ...
    use z; exact ...
    \end{leancode}
\item Provide an \lean{instance} to show that the completion of a densely normed field is still densely normed:
\begin{leancode}[hyperurl={https://github.com/YijunYuan/SphericalCompleteness/blob/717759d2b22a5720e74c14160a8385ac17f32abc/SphericalCompleteness/External/DenselyNormedField.lean\#L21-L38}]
instance instDenselyNormedFieldCompletion
    {α : Type*} [DenselyNormedField α] [CompletableTopField α] :
    DenselyNormedField (UniformSpace.Completion α)
\end{leancode}
\end{enumerate}
After these preparations, the instance
\lean{DenselyNormedField ℂ_[p]} can be automatically inferred by \lean{inferInstance}.
\subsubsection{Separability of $\bfC_p$}
It is easy to see, and is already formalized in Mathlib, that the completion of a separable metric space is still separable. Therefore, it suffices to show that $\bfQ^{\opn{alg}}_p$ is separable.

The approach is as follows:
\begin{enumerate}
    \item Take the algebraic closure $\bfQ^{\opn{alg}}$ of $\bfQ$ as a subset of $\bfQ^{\opn{alg}}_p$. Show that $\bfQ^{\opn{alg}}$ is countable. 
    
    Although it is mathematically trivial that $\bfQ^{\opn{alg}}$ is a subfield of $\bfQ^{\opn{alg}}_p$, this is not the case in Mathlib, as \lean{AlgebraicClosure ℚ} and \lean{AlgebraicClosure ℚ_[p]} are defined as two individual types. At first glance, one may want to use the fact that for any fields $F_1,F_2$, a non-zero morphism $F_1\lto F^{\opn{alg}}_2$ can always be extended to a morphism $F^{\opn{alg}}_1\lto F^{\opn{alg}}_2$. In particular, there exists an embedding $\bfQ^{\opn{alg}}\hookrightarrow\bfQ^{\opn{alg}}_p$, and one can then work with the image of this embedding.
    \begin{leancode}
Set.range <| @IsAlgClosed.lift (PadicAlgCl p) _ _ ℚ _ (AlgebraicClosure ℚ) _ _ _ _ _ _ _ -- : Set (PadicAlgCl p)
    \end{leancode}
In practice, we find the set $\left\{z\in \bfQ^{\opn{alg}}_p \vert z \text{ is algebraic over } \bfQ\right\}$
    \begin{leancode}[]
{z : PadicAlgCl p | IsAlgebraic ℚ z} 
\end{leancode}
\noindent more convenient to work with, as one can easily talk about the minimal polynomial of an element of this set over $\bfQ$, as well as its roots, in Lean without passing through any field morphisms. Moreover, the countability of this set can be proved in one line by \href{https://github.com/leanprover-community/mathlib4/blob/fabf563a7c95a166b8d7b6efca11c8b4dc9d911f/Mathlib/Algebra/AlgebraicCard.lean#L70-L74}{\lean{Algebraic.countable}} in Mathlib, which states that the set of algebraic elements (in a given extension field) over a countable field is countable.
    \item Fix $z\in\bfQ^{\opn{alg}}_p$ with minimal polynomial $f(X)\in\bfQ_p[X]$ over $\bfQ_p$. The \textbf{continuity of roots} (cf. \cite[§ 3.4.1, Proposition 1]{MR746961}) ensures that if we approximate $f$ closely enough by polynomials over $\bfQ$ with respect to the Gauss norm\footnote{For any polynomial $f=\sum_{i=0}^n a_i X^i\in\bfQ_p[X]$, the Gauss norm of $f$ is defined as $\max_{0 \leq i \leq n} \lVert a_i\rVert$.} (which is always possible since $\bfQ$ is dense in $\bfQ_p$), then $z$ can be approximated by roots of these polynomials, which are elements of $\bfQ^{\opn{alg}}$. This shows that $\bfQ^{\opn{alg}}$ is dense in $\bfQ^{\opn{alg}}_p$.
\end{enumerate}

\begin{remark}
    In the initial version of this project, the continuity of roots was implemented as an internal result \textsuperscript{\href{https://github.com/YijunYuan/SphericalCompleteness/blob/1d6181509701709e231d979cd2bfb7174d22f3e6/SphericalCompleteness/External/ContinuityOfRoots.lean\#L298-L460}{\faExternalLink*}} by the author. A more general variant of this result is now available in Mathlib, due to the work of Jiedong Jiang (cf. \cite{jjdcpalgclosed}), under the name \href{https://github.com/leanprover-community/mathlib4/blob/fabf563a7c95a166b8d7b6efca11c8b4dc9d911f/Mathlib/Analysis/Normed/Field/Approximation.lean#L130-L162}{\lean{exists_monic_and_natDegree_eq_and_norm_map_algebraMap_coeff_sub_lt}}. We then adopt this new version in our project to simplify the code.
\end{remark}

\section{Birkhoff–James Orthogonality}\label{sec:9468}
In this and the following sections, we once and for all fix a non-trivially normed field $\bbK$, two normed spaces $E,F$ over $\bbK$, and a $\bbK$-linear subspace $D$ of $E$.
\begin{leancode}
(𝕜 : Type*) [NontriviallyNormedField 𝕜]
(E : Type*) [SeminormedAddCommGroup E] [NormedSpace 𝕜 E] [IsUltrametricDist E]
(F : Type*) [SeminormedAddCommGroup F] [NormedSpace 𝕜 F] [IsUltrametricDist F]
(D : Submodule 𝕜 E)
\end{leancode}
\noindent All code snippets in this and the following sections are formalized under these assumptions, except for:
\begin{enumerate}
    \item some of the variables \lean{𝕜}, \lean{E}, \lean{F}, \lean{D} are made implicit to simplify the call site;
    \item in some places, one may need to replace \lean{SeminormedAddCommGroup} with \lean{NormedAddCommGroup} to ensure that the norm is non-degenerate, which is required by certain results.
\end{enumerate}
These differences will not be explicitly indicated in this article, for simplicity.

\subsection{Definition and basic properties}
In Hilbert spaces over $\bfR$ or $\bfC$, the orthogonality between two vectors $x,y$ is defined via the inner product. However, there is no inner product in the ultrametric world:
the parallelogram law fails in general\footnote{To see this, let $\norm{x}=\norm{y}=r>0$; then the ultrametric inequality implies $\norm{x+y}\leq r$ and $\norm{x-y}\leq r$. However, the parallelogram law requires that $\norm{x+y}^2+\norm{x-y}^2=4r^2$.}. The concept of \textbf{Birkhoff–James orthogonality} provides a way to define an orthogonality in normed spaces without inner products:
\begin{definition}[{cf. \cite{jamesOrthogonalityLinearFunctionals1947}}]
    Let $E$ be a normed space over a normed field $K$ (not necessarily ultrametric). For any $x,y\in E$, we say that $x$ is \textbf{Birkhoff–James orthogonal} to $y$, if $\lVert x\rVert\leq \lVert x+\lambda y\rVert$ is satisfied for every $\lambda\in K$.
\end{definition}
\begin{remark}
    When there is an inner product on $E$, the Birkhoff–James orthogonality coincides with the usual orthogonality defined via the inner product.
\end{remark}
In \cite{rooijNonArchimedeanFunctionalAnalysis1978} and in our project, the Birkhoff–James orthogonality is restated as follows:
\begin{definition}
    Let $E$ be a normed space over $\bbK$. Let $x,y$ be two elements in $E$ and let $V, W$ be two linear subspaces of $E$. We say that
    \begin{enumerate}
        \item $x$ is orthogonal to $V$, denoted by $x \perp_{\rmm} V$, if $\norm{x}=d(x,V)$;
        \item $x$ is orthogonal to $y$, denoted by $x \perp_{\rmm} y$, if $x \perp_{\rmm} [y]$, where $[y]$ is the linear subspace generated by $y$;
        \item $V$ is orthogonal to $W$, denoted by $V \perp_{\rms} W$, if for any $v\in V$, one has $v \perp_{\rmm} W$.
    \end{enumerate}
\end{definition}
This is formalized in our project as follows:
\begin{leancode}[hyperurl={https://github.com/YijunYuan/SphericalCompleteness/blob/717759d2b22a5720e74c14160a8385ac17f32abc/SphericalCompleteness/NormedVectorSpace/Orthogonal/Defs.lean}]
def IsMOrtho (x : E) (D : Subspace 𝕜 E) := Metric.infDist x D = ‖x‖

def IsVOrtho (x y : E) := Metric.infDist x (𝕜 ∙ y) = ‖x‖

def IsOrtho (F1 : Subspace 𝕜 E) (F2 : Subspace 𝕜 E) := ∀ x ∈ F1, IsMOrtho 𝕜 x F2

notation:50 x " ⟂ₘ " F => IsMOrtho _ x F
notation:50 F " ⟂ₛ " G => IsOrtho _ F G
notation:50 x " ⟂[" 𝕜 "] " y => IsVOrtho 𝕜 x y
\end{leancode}

Various properties of the Birkhoff–James orthogonality --- symmetry, homogeneity, etc. --- are formalized in our project. Among them, we prove that this formulation of orthogonality is equivalent to James's original definition:
\begin{leancode}[hyperurl={https://github.com/YijunYuan/SphericalCompleteness/blob/717759d2b22a5720e74c14160a8385ac17f32abc/SphericalCompleteness/NormedVectorSpace/Orthogonal/Basic.lean\#L83-L98}]
lemma IsVOrtho.iff_birkhoffJames (x y : E) : (x ⟂[𝕜] y) ↔ ∀ c : 𝕜, ‖x‖ ≤ ‖x + c • y‖
\end{leancode}

\subsection{Orthogonal projection and orthogonal complement}
\begin{definition}
    Let $E$ be a normed space over $\bbK$, and let $V$ be a subspace of $E$. 
    \begin{enumerate}
        \item A subspace $W$ of $E$ is called an \textbf{orthogonal complement} of $V$, if \begin{enumerate}
        \item $W$ is a complement of $V$, i.e. $V+W= E$ and $V\cap W=\{0\}$;
        \item $V \perp_{\rms} W$.
    \end{enumerate}
    \item Suppose $V$ admits an orthogonal complement $W$. The \textbf{orthogonal projection} onto $V$ along $W$ is the linear map $P\colon E\lto V$ satisfying $P\circ P = P$ and $\ker(P)=W$.
    \end{enumerate}
\end{definition}
As is indicated in the definition, the orthogonal complement may not exist in general. However, if the subspace is spherically complete, then one can always find an orthogonal complement:
\begin{theorem}[{cf. \cite[Corollary 4.7]{rooijNonArchimedeanFunctionalAnalysis1978}}]\label{thm:52486}
    Let $E$ be a normed space over $\bbK$, and let $D$ be a subspace of $E$. If $D$ is spherically complete, then there exists an orthogonal complement $D^{\perp}$ of $D$.
\end{theorem}
This theorem is a direct consequence of the following result:
\begin{proposition}\label{lem:6336}
    Let $E$ be a normed space over $\bbK$, and let $D$ be a subspace of $E$. If $D$ is spherically complete, then there exists a continuous linear map $T\colon E\to D$ such that $D$ is invariant under $T$ and $\norm{T}\leq 1$.
\end{proposition}
\begin{leancode}[hyperurl={https://github.com/YijunYuan/SphericalCompleteness/blob/717759d2b22a5720e74c14160a8385ac17f32abc/SphericalCompleteness/NormedVectorSpace/Orthogonal/OrthComp.lean\#L64-L87}]
theorem exists_orthProj_of_sphericallyComplete [SphericallyCompleteSpace D] :
    ∃ T : E →L[𝕜] ↥D, (∀ a ∈ D, T a = a) ∧ ‖T‖ ≤ 1
\end{leancode}
\begin{proof}
    This is immediate if one takes $\calU=\{0\}$ and $\epsilon_0=1$ in \Cref{prop:25689} -- a result we will explain in detail in the next section.
\end{proof}
\begin{proof}[Proof of \Cref{thm:52486}]
    Let $D^{\perp}$ be the kernel of $T$ in \Cref{lem:6336}. We claim that $D^{\perp}$ is an orthogonal complement of $D$ and $T$ is the orthogonal projection onto $D$ along $D^{\perp}$.

    If $x\in D\cap D^{\perp}$, then $x=T x =0$. On the other hand, for any $x\in E$, one has $x-T x\in D^{\perp}$ and $T x\in D$, thus $x=(x-T x)+T x\in D+D^{\perp}$. This shows that $D^{\perp}$ is a complement of $D$.

    To show that $D \perp_{\rms} D^{\perp}$, let $x\in D^{\perp}$ and $y\in D$. Then one has
    $$\norm{y}=\norm{T(x-y)}\leq\norm{x-y},$$
    which implies by the ultrametric inequality that $\norm{x}\leq \norm{x-y}$. As a result, one obtains $\norm{x}\leq \opn{dist}(x,D)$ and consequently $x\perp_{\rmm} D$.
\end{proof}

It is a common strategy of defining the data required in a proof, followed by the specifications of the said data in Lean formalizations. Specialized to the orthogonal projection and orthogonal complement, we first prove the existence of the continuous linear map $T$ in \Cref{lem:6336}, and then define the orthogonal projection and orthogonal complement via $T$.
\begin{leancode}[hyperurl={https://github.com/YijunYuan/SphericalCompleteness/blob/717759d2b22a5720e74c14160a8385ac17f32abc/SphericalCompleteness/NormedVectorSpace/Orthogonal/OrthComp.lean\#L180-L198}]
-- The orthogonal projection
def orthProj [SphericallyCompleteSpace D] : E →L[𝕜] ↥D :=
  (exists_orthProj_of_sphericallyComplete 𝕜 D).choose -- the map $T$ in $\text{\Cref{lem:6336}}$
\end{leancode}
\begin{leancode}[hyperurl={https://github.com/YijunYuan/SphericalCompleteness/blob/717759d2b22a5720e74c14160a8385ac17f32abc/SphericalCompleteness/NormedVectorSpace/Orthogonal/OrthComp.lean\#L89-L103}]
-- The orthogonal complement
def orthComp : Submodule 𝕜 E :=(exists_orthProj_of_sphericallyComplete 𝕜 D).choose
\end{leancode}
\noindent and formalize their properties as \lean{theorem}s:
\begin{leancode}[hyperurl={https://github.com/YijunYuan/SphericalCompleteness/blob/717759d2b22a5720e74c14160a8385ac17f32abc/SphericalCompleteness/NormedVectorSpace/Orthogonal/OrthComp.lean\#L126-L159}]
theorem IsOrtho.orthComp : (D ⟂ₛ (orthComp 𝕜 D))
\end{leancode}
\begin{leancode}[hyperurl={https://github.com/YijunYuan/SphericalCompleteness/blob/717759d2b22a5720e74c14160a8385ac17f32abc/SphericalCompleteness/NormedVectorSpace/Orthogonal/OrthComp.lean\#L105-L122}]
theorem isCompl_orthComp : IsCompl D (orthComp 𝕜 D)
\end{leancode}
\begin{leancode}[hyperurl={https://github.com/YijunYuan/SphericalCompleteness/blob/717759d2b22a5720e74c14160a8385ac17f32abc/SphericalCompleteness/NormedVectorSpace/Orthogonal/OrthComp.lean\#L200-L222}]
theorem norm_orthProj_le_one : ‖orthProj 𝕜 D‖ ≤ 1

theorem orthProj_id : ∀ a ∈ D, (orthProj 𝕜 D) a = a
\end{leancode}

\subsection{Application: Finite dimensional normed spaces}
As an application — seemingly unrelated to the orthogonality at first glance — we prove the following result:
\begin{theorem}
    If $\bbK$ is spherically complete, then every finite dimensional normed space over $\bbK$ is spherically complete.
\end{theorem}
\begin{leancode}[hyperurl={https://github.com/YijunYuan/SphericalCompleteness/blob/717759d2b22a5720e74c14160a8385ac17f32abc/SphericalCompleteness/NormedVectorSpace/Basic.lean\#L136-L164}]
theorem of_finiteDimensional [FiniteDimensional 𝕜 E] : SphericallyCompleteSpace E
\end{leancode}
The strategy is a direct induction on the dimension:
\begin{enumerate}
    \item The $0$-dimensional case is trivial.
    \item The $1$-dimensional spaces are isometric to the base field, which is spherically complete by assumption.
    \begin{leancode}[hyperurl={https://github.com/YijunYuan/SphericalCompleteness/blob/717759d2b22a5720e74c14160a8385ac17f32abc/SphericalCompleteness/NormedVectorSpace/Basic.lean\#L49-L99}]
instance instSphericallyCompleteSpaceSpanSingleton
    [SphericallyCompleteSpace 𝕜] (z : E) :
    SphericallyCompleteSpace (Submodule.span 𝕜 {z}) 
    \end{leancode}
    \item Assume there exists an $n$-dimensional spherically complete subspace $D$ of $E$, where $n <\opn{dim}E$. Since $D$ is a proper subspace, one may pick, by \Cref{prop:3029}, an element $a\in E$ with $a\notin D$ and $a\perp_{\rms} D$.
    \item Since $D$ and $[a]$ are both spherically complete, so is their product space $D\times [a]$ (cf. \Cref{prop:1411}). The orthogonality $a\perp_{\rms} D$ makes the isomorphism
    \begin{equation}\label{eq:23053}
        D\oplus[a]\cong D\times [a]
    \end{equation}
    an isometry (again by \Cref{prop:3029}), so that spherical completeness transfers from $D\times[a]$ to $D\oplus[a]$ via (2) of \Cref{prop:1411}.
\end{enumerate}
Without $a\perp_{\rms} D$ the isomorphism \eqref{eq:23053} need not be an isometry, and spherical completeness would not transfer. The existence of such an $a$ and the resulting isometry are exactly the following proposition:
\begin{proposition}[{cf. \cite[page 57 \& Lemma 3.14]{rooijNonArchimedeanFunctionalAnalysis1978}}]\label{prop:3029}
    Let $E$ be a normed space over $\bbK$, and let $D$ be a spherically complete proper subspace of $E$. 
    \begin{enumerate}
        \item There exists an element $a_0\in E$ such that $a_0\notin D$ and $a_0\perp D$.
        \item One has the isometric isomorphism $D\times [a_0] \cong D+[a_0]$.
    \end{enumerate}
\end{proposition}
\begin{leancode}[hyperurl={https://github.com/YijunYuan/SphericalCompleteness/blob/717759d2b22a5720e74c14160a8385ac17f32abc/SphericalCompleteness/NormedVectorSpace/Orthogonal/IsMOrtho.lean\#L162-L254}]
theorem exists_isMOrtho_vec_of_finrank_lt
    [SphericallyCompleteSpace D] (hF : Module.finrank 𝕜 D < Module.finrank 𝕜 E) :
    ∃ (x : E), x ≠ 0 ∧ (x ⟂ₘ D)
\end{leancode}
\begin{leancode}[hyperurl={https://github.com/YijunYuan/SphericalCompleteness/blob/717759d2b22a5720e74c14160a8385ac17f32abc/SphericalCompleteness/NormedVectorSpace/Orthogonal/IsMOrtho.lean\#L40-L125}]
def IsMOrtho.directProdIsoSum (x : E) (hxF : x ⟂ₘ D) :
    (Submodule.span 𝕜 {x}) × D≃ₛₗᵢ[RingHom.id 𝕜] (Submodule.span 𝕜 {x}) + D
\end{leancode}

The isometric isomorphism $D\times[a_0]\cong D+[a_0]$ in \Cref{prop:3029} is not free: it hinges on the spherical completeness of $D$. This is easy to overlook in informal writing, where a single symbol $\cong$ is used for isomorphisms living in different categories. In the category of vector spaces, a direct sum is always naturally isomorphic to the direct product; but for normed spaces the ``same'' isomorphism need not be an isometry. Here formalization earns its keep: a strict type system forbids conflating the two $\cong$'s and forces one to pin down the category of each.

\section{Ultrametric Hahn–Banach theorem}\label{sec:45266}
In functional analysis, the Hahn–Banach theorem is a fundamental result about the extension of bounded linear functionals. There are various versions of this theorem in different contexts. One of them is stated as follows:
\begin{theorem}
    Every continuous linear functional defined on a vector subspace $M$ of a normed space $X$ over $\bfR$ or $\bfC$ has a continuous linear extension $F$ to all of $X$ that satisfies $\lVert f\rVert=\lVert F\rVert$.
\end{theorem}

This is already formalized in Mathlib:
\begin{leancode}[hyperurl={https://github.com/leanprover-community/mathlib4/blob/fabf563a7c95a166b8d7b6efca11c8b4dc9d911f/Mathlib/Analysis/Normed/Module/HahnBanach.lean\#L43-L53}]
theorem exists_extension_norm_eq
    {𝕜 : Type u_1} [NontriviallyNormedField 𝕜] [IsRCLikeNormedField 𝕜]
    {E : Type u_2} [SeminormedAddCommGroup E] [NormedSpace 𝕜 E]
    (p : Subspace 𝕜 E) (f : StrongDual 𝕜 ↥p) :
    ∃ (g : StrongDual 𝕜 E), (∀ (x : ↥p), g ↑x = f x) ∧ ‖g‖ = ‖f‖    
\end{leancode}
\begin{remark}
    In the above theorem, \lean{IsRCLikeNormedField} is the typeclass in Mathlib saying that the norm field structure is the same as $\bfR$ or $\bfC$, and \lean{StrongDual} is the space of continuous linear functionals on $M$. 
\end{remark}

Unfortunately, this result does not hold if one simply replaces the base field $\bfR$ or $\bfC$ by an ultrametric normed field. Instead, one has to take spherical completeness into account:
\begin{theorem}[{cf. \cite[Theorem 4.8]{rooijNonArchimedeanFunctionalAnalysis1978}}]\label{thm:14044}
    Let $E,F$ be two normed spaces over $\bbK$, and let $D$ be a subspace of $E$. If $D$ is spherically complete or $F$ is spherically complete, then every $f\in\calL(D,F)$ can be extended to some $F\in\calL(E,F)$ such that $\lVert F\rVert=\lVert f\rVert$.
\end{theorem}
Instead of formalizing this theorem as two individual \lean{theorem}s in Lean, corresponding to the conditions \lean{SphericallyCompleteSpace D} and \lean{SphericallyCompleteSpace F} respectively, we introduce the following concept:
\begin{definition}
    The data $(D,F)$ is called a \textbf{Hahn–Banach pair} if for every $f\in\calL(D,F)$, there exists an extension $F\in\calL(E,F)$ of $f$ such that $\lVert F\rVert=\lVert f\rVert$.
    \begin{leancode}[hyperurl={https://github.com/YijunYuan/SphericalCompleteness/blob/717759d2b22a5720e74c14160a8385ac17f32abc/SphericalCompleteness/NormedVectorSpace/ContinuousLinearMap/HahnBanach.lean\#L56-L65}]
class IsHahnBanachPair [IsUltrametricDist E] [IsUltrametricDist F] : Prop where
  is_hb : ∀ f : D →L[𝕜] F, ∃ f' : E →L[𝕜] F, 
    (∀ v : E, (hv : v ∈ D) → f' v = f ⟨v, hv⟩) ∧ ‖f'‖ = ‖f‖
\end{leancode}
\end{definition}

We provide two \lean{instance}s to show that \lean{IsHahnBanachPair} holds when either $D$ or $F$ is spherically complete, and this is where the proof of \Cref{thm:14044} is actually formalized in Lean:
\begin{leancode}[hyperurl={https://github.com/YijunYuan/SphericalCompleteness/blob/717759d2b22a5720e74c14160a8385ac17f32abc/SphericalCompleteness/NormedVectorSpace/ContinuousLinearMap/HahnBanach.lean\#L77-L120}]
instance instIsHahnBanachPairOfDomainSphericalComplete [SphericallyCompleteSpace D] : 
    IsHahnBanachPair D F

instance instIsHahnBanachPairOfCodomainSphericalComplete
    [SphericallyCompleteSpace F] : IsHahnBanachPair D F
\end{leancode}
\noindent The statement of \Cref{thm:14044} then immediately follows from the definition of \lean{IsHahnBanachPair}.
\begin{leancode}[hyperurl={https://github.com/YijunYuan/SphericalCompleteness/blob/717759d2b22a5720e74c14160a8385ac17f32abc/SphericalCompleteness/NormedVectorSpace/ContinuousLinearMap/HahnBanach.lean\#L67-L75}]
theorem hahn_banach [IsHahnBanachPair D F] (f : D →L[𝕜] F) :
    ∃ f' : E →L[𝕜] F, (∀ v : E, (hv : v ∈ D) → f' v = f ⟨v, hv⟩) ∧ ‖f'‖ = ‖f‖ :=
  IsHahnBanachPair.is_hb f
\end{leancode}
\begin{remark}
    A normed space $F$ over $\bbK$ is called \textbf{injective} if $(D,F)$ is a Hahn–Banach pair for any normed space $E$ and any subspace $D$ of $E$. Ingleton proved in 1952 that a normed space $F$ is injective if and only if $F$ is spherically complete (cf. \cite[Theorem 2]{ingletonHahnBanachTheoremNonArchimedeanvalued1952}), which implies \Cref{thm:14044} in the case of a spherically complete codomain $F$. For the case of a spherically complete domain $D$, we were unable to locate the first appearance of this result in the literature. The author learned this result from van Rooij's book \cite[Theorem 4.8]{rooijNonArchimedeanFunctionalAnalysis1978}.
\end{remark}

%In the rest of this section, we explain the formalization of the proof of \Cref{thm:14044}.

%\subsection{Formalization of a key lemma}
The proof of \Cref{thm:14044} heavily relies on the following key result:
\begin{proposition}[{cf. \cite[Lemma 4.4]{rooijNonArchimedeanFunctionalAnalysis1978}}]\label{prop:25689}
    Let $E,F$ be two normed spaces over $\bbK$, with $F$ spherically complete. Let $D$ be a subspace of $E$, $S\in\calL(D,F)$. Let $\calU$ be a nonempty subset of $\calL(E,F)$ and for each $U\in\calU$, let $\epsilon_U>0$ be so that
    $$\lVert U-V\rVert\leq \max\{\epsilon_U,\epsilon_V\}, \text{ for every } U, V\in\calU,$$
    $$\lVert Sx-Ux\rVert\leq \epsilon_U\lVert x\rVert,\text{ for every } U\in\calU \text{ and } x\in D.$$
    Then $S$ has an extension $\overline{S}\in\calL(E,F)$ such that $\lVert \overline{S}-U\rVert\leq \epsilon_U$ for every $U\in\calU$.
\end{proposition}
\begin{leancode}[hyperurl={https://github.com/YijunYuan/SphericalCompleteness/blob/717759d2b22a5720e74c14160a8385ac17f32abc/SphericalCompleteness/NormedVectorSpace/ContinuousLinearMap/Extension.lean\#L527-L608}]
theorem exists_extension_opNorm_le (S : D →L[𝕜] F)
    {𝒰 : Set (E →L[𝕜] F)} (h𝒰 : 𝒰.Nonempty)
    (ε : ↑𝒰 → ℝ) (hε1 : ∀ T : ↑𝒰, 0 < ε T)
    (hε2 : ∀ U V : ↑𝒰, ‖U.val - V.val‖ ≤ max (ε U) (ε V))
    (hε3 : ∀ U : ↑𝒰, ∀ x : D, ‖S x - U.val x‖ ≤ ε U * ‖x‖) :
    ∃ (T : E →L[𝕜] F), (∀ x : D, T x = S x) ∧ (∀ U : ↑𝒰, ‖T - U.val‖ ≤ ε U)
\end{leancode}

\begin{proof}[Proof of \Cref{thm:14044}]
    When the codomain $F$ is spherically complete, take $\calU=\{0\}$ and $\epsilon_0=\lVert S\rVert$ in \Cref{prop:25689} and the resulting map $\overline{S}$ is the desired extension. 

    When $D$ is spherically complete, the composition $S\circ T$ gives the extension, where $T$ is the orthogonal projection onto $D$ along $D^\perp$.
\end{proof}
As this proof is simple, we skip the formalization details here and focus on the formalization of \Cref{prop:25689} in the following subsections. In particular, in the rest of this section, the meaning of variables \lean{S}, \lean{𝒰}, \lean{ε}, \lean{hε1}, \lean{hε2}, \lean{hε3} will keep the same as in \lean{exists_extension_opNorm_le}, and their declaration will be omitted for simplicity.

\subsection{Reduce to the case of codimension 1}
In the proof of \cite[Lemma 4.4]{rooijNonArchimedeanFunctionalAnalysis1978}, van Rooij writes:
\begin{quote}
    By a standard application of Zorn's lemma we may assume that there exists an $a\in E$, $a\notin D$ such that $E=D+[a]$.
\end{quote}
To translate this assertion into a formalized proof, we need to clarify the partially ordered set $\frakP$ we are working with. In this context, $\frakP$ consists of the pairs $(M,T)$, where $M$ is a submodule of $E$ and $T\in\calL(M,K)$ is a continuous linear map, satisfying the following conditions:
\begin{enumerate}
    \item for any $x\in D$, $T x = S x$;
    \item for any $x\in M$ and $U\in\calU$, $\lVert T x - U x\rVert\leq \epsilon_U\lVert x\rVert$;
    \item $(M_1,T_1)\leq (M_2,T_2)$ if $M_1\subseteq M_2$ and $T_2$ extends $T_1$.
\end{enumerate}
In our project, this set $\frakP$ is formalized as a \lean{structure}:
\begin{leancode}[hyperurl={https://github.com/YijunYuan/SphericalCompleteness/blob/717759d2b22a5720e74c14160a8385ac17f32abc/SphericalCompleteness/NormedVectorSpace/ContinuousLinearMap/Extension.lean\#L328-L351}]
structure PartialExtension where
  M : Submodule 𝕜 E
  hDM : D ≤ M -- D is a submodule of M
  T : M →L[𝕜] F
  hT : ∀ x : D, T ⟨x, hDM x.prop⟩ = S x -- Condition 1
  hU : ∀ U : ↑𝒰, ∀ x : M, ‖T x- U.val x‖ ≤ (ε U) * ‖x‖ -- Condition 2
\end{leancode}
\noindent And we formalize the partial order on $\frakP$ as an \lean{instance}:
\begin{leancode}[hyperurl={https://github.com/YijunYuan/SphericalCompleteness/blob/717759d2b22a5720e74c14160a8385ac17f32abc/SphericalCompleteness/NormedVectorSpace/ContinuousLinearMap/Extension.lean\#L370-L383}]
instance instPartialOrderPartialExtension : -- Condition 3
    PartialOrder (PartialExtension 𝕜 E F S 𝒰 h𝒰 ε) where
    le a b := ∃ hab : a.M ≤ b.M , ∀ x : a.M, b.T ⟨x.val, hab x.prop⟩ = a.T x    
    ... -- Proof of reflexivity, antisymmetry and transitivity. Omitted.
\end{leancode}

To apply Zorn's lemma, we need to show that every chain $\{(M_\alpha,T_\alpha)\}_{\alpha\in I} \subseteq \frakP$ has an upper bound. It is easy to verify that if one takes $M_{\max}$ to be the union of all the $M_\alpha$'s, then $M_{\max}$ is still a submodule of $E$, and the $T_\alpha$'s can be glued together to form a linear map $T_{\max}\in\calL(M_{\max},F)$. Moreover, we have:
\begin{lemma}
    The linear map $T_{\max}$ is continuous. In particular, $(M_{\max},T_{\max})$ is an upper bound of the chain $\{(M_\alpha,T_\alpha)\}_{\alpha\in I}$.
\end{lemma}
\begin{proof}
    Take an arbitrary $U\in\calU$ and set $C\coloneqq \max\{\epsilon_U,\lVert U\rVert\}$. For any $x\in M_{\max}$, there exists some $\alpha\in I$ such that $x\in M_\alpha$. One computes that
    $$\norm{T_{\max}x}=\norm{T_\alpha x}\leq \max\{\norm{T_\alpha x - U x}, \norm{U x}\}\leq \max\{\epsilon_U\norm{x},\norm{U}\norm{x}\}\leq C \norm{x}.$$
    This shows that $T_{\max}$ is bounded, hence continuous.
\end{proof}
When formalizing the fact that $T_{\max}$ is linear and continuous, the locally-defined nature of $T_{\max}$ forces us to repeatedly use \lean{Classical.choose_spec} and the fact that chains are directed in the formalized proof, which makes the code quite verbose and less transparent. Notwithstanding these challenges, we successfully formalize the existence of upper bounds in our project:
\begin{leancode}[hyperurl={https://github.com/YijunYuan/SphericalCompleteness/blob/717759d2b22a5720e74c14160a8385ac17f32abc/SphericalCompleteness/NormedVectorSpace/ContinuousLinearMap/Extension.lean\#L485-L522}]
lemma bddAbove_of_chain_of_partial_extension
    (P : Set (PartialExtension 𝕜 E F S 𝒰 h𝒰 ε))
    (hP : IsChain (fun x1 x2 ↦ x1 ≤ x2) P) (hhP : P.Nonempty) : BddAbove P
\end{leancode}

By applying Zorn's lemma, we obtain a maximal element $(M_{\max},T_{\max})$ in $\frakP$. This reduces the proof of \Cref{thm:14044} to the case where $E$ is of codimension 1 over $D$:
\begin{lemma}
    If \Cref{thm:14044} holds when $E=D+[a]$ for some $a\in E\backslash D$, then it holds in general.
\end{lemma}
\begin{proof}
    We claim that $M_{\max}=E$. If not, take $a\in E\backslash M_{\max}$, and set $M'\coloneqq M_{\max}+[a]$. By assumption, there exists some $T'\in\calL(M',F)$ extending $T_{\max}$ such that for any $U\in\calU$ and $x\in M'$, $\lVert T'x - U x\rVert\leq \epsilon_U\lVert x\rVert$. This shows that $(M',T')\in\frakP$ and $(M_{\max},T_{\max})<(M',T')$, contradicting the maximality of $(M_{\max},T_{\max})$.
\end{proof}
The corresponding formalized proof of this lemma is integrated in the proof of \Cref{prop:25689}, which we omit here for brevity.

\subsection{Remarks on the codimension 1 case}
Now we assume that there exists some $a\in E\backslash D$ such that $E=D+[a]$.

It is shown in \cite{rooijNonArchimedeanFunctionalAnalysis1978} that the spherical completeness of $F$ implies the existence of an element $z_0\in F$ such that 
$$\norm{S x+z_0 - U(x + a)}\leq \epsilon_U\norm{x + a}$$
holds for every $U\in\calU$ and $x\in D$.
\begin{leancode}[hyperurl={https://github.com/YijunYuan/SphericalCompleteness/blob/717759d2b22a5720e74c14160a8385ac17f32abc/SphericalCompleteness/NormedVectorSpace/ContinuousLinearMap/Extension.lean\#L96-L186}]
lemma exists_extensionValue_norm_le [SphericallyCompleteSpace F] 
    {a : E} (ha1 : a ∉ D) :
    ∃ z0 : F, ∀ (x : ↥D) (U : ↑𝒰), ‖S x + z0 - U.val (↑x + a)‖ ≤ ε U * ‖↑x + a‖
\end{leancode}
\noindent Then the desired extension is given by
$$\overline{S}\colon x + \lambda\cdot a \longmapsto S x + \lambda\cdot z_0.$$
\begin{leancode}[hyperurl={https://github.com/YijunYuan/SphericalCompleteness/blob/717759d2b22a5720e74c14160a8385ac17f32abc/SphericalCompleteness/NormedVectorSpace/ContinuousLinearMap/Extension.lean\#L205-L216}]
def codimOneExtension [SphericallyCompleteSpace F] {a : E} (ha1 : a ∉ D) :
    (D + Submodule.span 𝕜 {a}) → F := fun M => by
    have := Submodule.mem_sup.1 M.prop -- $M = x + y$ with $x \in D$ and $y\in [a]$
    let lambda := (Submodule.mem_span_singleton.1 
      this.choose_spec.2.choose_spec.1).choose -- There exists $\lambda\in\ $𝕜 such that $y = \lambda\cdot a$
    use S ⟨this.choose, this.choose_spec.1⟩ +
      lambda • (exists_extensionValue_norm_le ha1 S h𝒰 hε1 hε2 hε3).choose
\end{leancode}
\begin{remark}
    In the above code snippet, for any $M\in D+[a]$, the corresponding $x\in D$ and $\lambda\in K$ are obtained by applying \lean{Classical.choose_spec} and \lean{Classical.choose} in a mixed manner. This makes it hard to reason about the properties of this map. One potential remedy is to introduce local notations for the required data during the proof. This makes the code more readable, at the cost of a longer preamble (which is not a problem on most occasions).
\end{remark}
We prove that $\overline{S}$ is continuous and linear by showing that it is a bounded linear map:
\begin{leancode}[hyperurl={https://github.com/YijunYuan/SphericalCompleteness/blob/717759d2b22a5720e74c14160a8385ac17f32abc/SphericalCompleteness/NormedVectorSpace/ContinuousLinearMap/Extension.lean\#L256-L281}]
def isBoundedLinearMap_codimOneExtension [SphericallyCompleteSpace F]
    {a : E} (ha1 : a ∉ D) :
    IsBoundedLinearMap 𝕜 (codimOneExtension ha1 S h𝒰 hε1 hε2 hε3)
\end{leancode}
\noindent The structure \href{https://github.com/leanprover-community/mathlib4/blob/fabf563a7c95a166b8d7b6efca11c8b4dc9d911f/Mathlib/Topology/Algebra/Module/ContinuousLinearMap/Basic.lean#L58-L65}{\lean{ContinuousLinearMap}} can then be obtained by \href{https://github.com/leanprover-community/mathlib4/blob/fabf563a7c95a166b8d7b6efca11c8b4dc9d911f/Mathlib/Analysis/Normed/Operator/BoundedLinearMaps.lean#L102-L107}{\lean{IsBoundedLinearMap.toContinuousLinearMap}}.
\iffalse
, as one has to write code like
\begin{leancode}[hyperurl={https://github.com/YijunYuan/SphericalCompleteness/blob/1d6181509701709e231d979cd2bfb7174d22f3e6/SphericalCompleteness/NormedVectorSpace/ContinuousLinearMap/SupportingResults.lean\#L278-L280}]
have tt := (rooij_lemma_4_4_z0_prop ha1 S h𝒰 hε1 hε2 hε3)
      ⟨(Submodule.mem_sup.1 x.prop).choose, (Submodule.mem_sup.1 x.prop).choose_spec.1⟩
      ((Submodule.mem_span_singleton.1
      (Submodule.mem_sup.1 x.prop).choose_spec.2.choose_spec.1).choose) ⟨h𝒰.some, h𝒰.some_mem⟩    
\end{leancode}
\noindent which is hard to read and debug.

\begin{remark}
    At first glance, it seems that this predicament can be resolved by using \lean{rcases} tactic to extract $x$ and $\lambda$ from $M$, which is more elegant and readable. However, this is not the case, as \lean{rcases} tactic is not usable in the definition of general functions. As a result, if one invokes \lean{rcases} tactic to prove the properties of \lean{rooij_lemma_4_4_T}, then Lean is unable to unify these \lean{Classical.choose}-based terms and \lean{rcases}-resulted terms.
\end{remark}
\fi

\subsection{Application: Spherical completeness of $\calL(E,F)$}
Besides its use in constructing the orthogonal complement (cf. \Cref{lem:6336}) and proving the ultrametric Hahn–Banach theorem (\Cref{thm:14044}), \Cref{prop:25689} can also be applied to show the spherical completeness of the space of continuous linear maps to a spherically complete codomain:
\begin{theorem}[{cf. \cite[Corollary 4.5]{rooijNonArchimedeanFunctionalAnalysis1978}}]\label{thm:25444}
    Let $E,F$ be two ultrametric normed spaces over an ultrametric normed field $K$. If $F$ is spherically complete, then so is $\calL(E,F)$.
\end{theorem}
The proof of \Cref{thm:25444} in \cite{rooijNonArchimedeanFunctionalAnalysis1978} is concise:
\begin{quote}
    Let $B(T_1,\epsilon_1)\supset B(T_2,\epsilon_2)\supset\cdots$ be balls in $\calL(E,F)$. Apply \Cref{prop:25689}, taking $D=\{0\}$, $\calU=\{T_n\colon n\in\bfN\}$, $\epsilon_{T_n}=\epsilon_n$.
\end{quote}
During the formalization, we recognized that there are extra details that need to be taken care of:
\begin{enumerate}
    \item To apply \Cref{prop:25689}, one has to verify that $\epsilon_U>0$
    for every $U\in\calU$ in the current context. To do so, the extra condition that the radii of these balls are strictly decreasing is required. By \Cref{thm:5601}, this extra condition does not incur any loss of generality.
    \item The second point is more serious. The declaration ``$\epsilon_{T_n}=\epsilon_n$'' implicitly uses the composition
    $$\calU\lto \bfN\lto \bfR,\ T_n\longmapsto n\longmapsto \epsilon_n,$$
    which is not well-defined in general, as the map $n\longmapsto T_n$ may not be injective. To bypass this issue, we consider the condition:
    \begin{equation}\label{eq:58585}
        \text{\bfseries\itshape for any $n\in\bfN$, there exists $N\in\bfN$ such that for any $i>N$, $T_n\neq T_i$,}
    \end{equation}
    or, equivalently, any element of $\calU$ only appears finitely many times in the sequence $\{T_n\}_{n\in\bfN}$.

    If this condition holds, then one can take a subsequence $\{T_{n_k}\}_{k\in\bfN}$ of $\{T_n\}_{n\in\bfN}$ such that the map $k\longmapsto T_{n_k}$ is injective.
    \begin{leancode}[hyperurl={https://github.com/YijunYuan/SphericalCompleteness/blob/717759d2b22a5720e74c14160a8385ac17f32abc/SphericalCompleteness/External/Sequence.lean\#L82-L111}]
theorem exists_injective_subseq_of_finite_duplication {α : Type*} (seq : ℕ → α)
    (hseq : ∀ n : ℕ, ∃ N, ∀ i > N, seq n ≠ seq i) :
    ∃ φ : ℕ → ℕ, StrictMono φ ∧ Function.Injective (seq ∘ φ)
    \end{leancode}
    \noindent Then van Rooij's original argument can be applied to the balls $\{B(T_{n_k},\epsilon_{n_k})\}_{k\in\bfN}$, which yields a desired element in the intersection of all the $B(T_{n_k},\epsilon_{n_k})$ and consequently of all the $B(T_n,\epsilon_n)$.

    Conversely, if \eqref{eq:58585} does not hold, then there exists some $T\in\calU$ such that $T=T_n$ for infinitely many $n$. Then $T$ itself witnesses the non-emptiness of the intersection of all the $B(T_n,\epsilon_n)$.
\end{enumerate}

This kind of subtlety can easily be overlooked in traditional mathematical writing, but has to be addressed carefully in a formalized proof. In our project, we successfully formalize the proof of \Cref{thm:25444} as an \lean{instance} by taking these details into account.

\section{Spherical completion}\label{sec:27285}
The completion of a metric space can be characterized via the following nice universal property:
\begin{theorem}[{cf. \cite[Chapter II, §3.8, Proposition 16]{bourbakiGeneralTopologyChapters1995}}]\label{thm:14845}
    Let $(X,d)$ be a metric space. The completion of $X$ is a complete metric space $\widehat{X}$ together with an isometric embedding $i\colon X\to \widehat{X}$ such that for any complete metric space $Y$ and any uniformly continuous map $f\colon X\to Y$, there exists a unique uniformly continuous map $\widehat{f}\colon \widehat{X}\to Y$ such that $\widehat{f}\circ i=f$.
\end{theorem}
Moreover, the completion can be constructed in a very explicit way, e.g. as the space of all Cauchy sequences on $X$ modulo the relation of being ``infinitesimally close''.

Unfortunately, there is no canonical way to construct such a ``spherical completion'' for an arbitrary metric space. On the other hand, for ultrametric normed vector spaces, van Rooij formulated the following concept:
\begin{definition}[{cf. \cite[page 146]{rooijNonArchimedeanFunctionalAnalysis1978}}]\label{def:}
    Let $E$ be a normed space over $\bbK$. A \textbf{spherical completion} of $E$ is a spherically complete normed space $F$ over $\bbK$ together with an isometric linear embedding $i\colon E\to F$ such that $F$ has no proper spherically complete subspace containing $i(E)$.
\end{definition}
This is designed as a \lean{class} in our project:
\begin{leancode}[hyperurl={https://github.com/YijunYuan/SphericalCompleteness/blob/717759d2b22a5720e74c14160a8385ac17f32abc/SphericalCompleteness/NormedVectorSpace/SphericalCompletion/Defs.lean\#L43-L67}]
class IsSphericalCompletion (ι : E →ₗᵢ[𝕜] F) : Prop 
    extends SphericallyCompleteSpace F where
  minimal_of_sphericallyComplete :
    ∀ D : Submodule 𝕜 F, SphericallyCompleteSpace D → ι.range ≤ D → D = ⊤
\end{leancode}

The main goal of this section is to formalize the existence and uniqueness of spherical completions, and to establish several of their properties. The following result of Fleischer is of particular importance:
\begin{theorem}[{cf. \cite{fleischerMaximalityUltracompletenessNormed1958}, \cite[Theorem 4.43]{rooijNonArchimedeanFunctionalAnalysis1978}}]\label{thm:63198}
    Every normed vector space over an ultrametric normed field admits a spherical completion. Any two spherical completions of a normed vector space are isometrically isomorphic.
\end{theorem}
This is formalized as two separate \lean{theorem}s in our project, one for existence and the other for uniqueness:
\begin{leancode}[hyperurl={https://github.com/YijunYuan/SphericalCompleteness/blob/717759d2b22a5720e74c14160a8385ac17f32abc/SphericalCompleteness/NormedVectorSpace/SphericalCompletion/Basic.lean\#L109-L133}]
theorem exists_isSphericalCompletion : ∃ E₀ : Type u,
    ∃ _ : NormedAddCommGroup E₀, ∃ _ : NormedSpace 𝕜 E₀, ∃ _ : IsUltrametricDist E₀, 
    ∃ ι : E →ₗᵢ[𝕜] E₀, IsSphericalCompletion ι
\end{leancode}
\begin{leancode}[hyperurl={https://github.com/YijunYuan/SphericalCompleteness/blob/717759d2b22a5720e74c14160a8385ac17f32abc/SphericalCompleteness/NormedVectorSpace/SphericalCompletion/Basic.lean\#L189-L210}]
theorem IsSphericalCompletion.unique 
    {F₁ : Type*} [NormedAddCommGroup F₁] [NormedSpace 𝕜 F₁] [IsUltrametricDist F₁]
    {F₂ : Type*} [NormedAddCommGroup F₂] [NormedSpace 𝕜 F₂] [IsUltrametricDist F₂]
    {ι₁ : E →ₗᵢ[𝕜] F₁} {ι₂ : E →ₗᵢ[𝕜] F₂} :
    IsSphericalCompletion ι₁ → IsSphericalCompletion ι₂ →
    ∃ f : F₁ ≃ₗᵢ[𝕜] F₂, f.toLinearIsometry.comp ι₁ = ι₂
\end{leancode}
\begin{remark}
    Our formalized code proves the existence of the spherical completion in the same universe level \lean{u} as the original space $E$. This is a reasonable constraint, as the spherical completion is constructed as a subspace of a certain subquotient of a countable product of $E$ (cf. \Cref{prop:58585}), which is still at the same universe level as $E$.
\end{remark}
\begin{remark}
    The spherical completion does not admit as good a universal property as the usual completion of metric spaces: for two spherical completions $(\breve{E_1},i_1)$ and $(\breve{E_2},i_2)$ of $E$, the isometric isomorphism $\breve{E_1}\cong \breve{E_2}$ is not necessarily unique. See also \Cref{coro:23722}.
\end{remark}

\subsection{Immediate extensions}
The concept of immediate extensions is closely related to the spherical completeness of ultrametric normed spaces, and plays a key role in the proof of \Cref{thm:63198}. In \cite{rooijNonArchimedeanFunctionalAnalysis1978}, van Rooij defines immediate extensions as follows:
\begin{definition}\label{def:1145}
    Let $E$ be a normed space over $\bbK$ and let $D$ be a subspace of $E$. $E$ is called an \textbf{immediate extension} of $D$ if $0$ is the only element of $E$ that is orthogonal to $D$.
\end{definition}
As an experiment, we tried to formalize this definition using Claude Opus, a large language model developed by Anthropic. It gave a literal translation:
\begin{leancode}[hyperurl={https://github.com/YijunYuan/SphericalCompleteness/blob/717759d2b22a5720e74c14160a8385ac17f32abc/SphericalCompleteness/NormedVectorSpace/Immediate.lean\#L72-L75}]
def IsImmediateExtension (D : Submodule 𝕜 E) : Prop := ∀ v : E, (v ⟂ₘ D) → v = 0
\end{leancode}
While this formalization is mathematically correct, it is not a good design choice. One may need to handle two independent immediate extensions $E_1/D$, $E_2/D$ of $D$ in the proof of \Cref{thm:63198}. In this case, if one defines $D$ to be of type \lean{Submodule 𝕜 E₁}, then an awkward coercion is required to express the condition that $E_2$ is an immediate extension of $D$, which makes the code less concise. To resolve this issue, we reformulate \Cref{def:1145} in the following more flexible way:
\begin{definition}\label{def:8095}
    Let $E, F$ be normed spaces over $\bbK$. Let $f\colon E\lto F$ be an isometric linear map. We say $f$ is \textbf{immediate} if the only element of $F$ that is orthogonal to $\opn{Im}(f)$ is $0$.
\end{definition}
\begin{leancode}[hyperurl={https://github.com/YijunYuan/SphericalCompleteness/blob/717759d2b22a5720e74c14160a8385ac17f32abc/SphericalCompleteness/NormedVectorSpace/Immediate.lean\#L56-L70}]
def IsImmediate (f : E →ₗᵢ[𝕜] F) : Prop := ∀ v : F, (v ⟂ₘ LinearMap.range f) → v = 0
\end{leancode}
Compared to \Cref{def:1145}, this definition allows $E$ and $F$ to be independent types, which makes the statements of related results more convenient to work with. For example, van Rooij proved the following result, which is crucial for the uniqueness argument in \Cref{thm:63198}:
\begin{lemma}[{cf. \cite[Lemma 4.42]{rooijNonArchimedeanFunctionalAnalysis1978}}]
    Let $F$ be an immediate extension of a normed space $E$ and let $H$ be any spherically complete normed space containing $E$. Then $\opn{id}_E$ can be extended to a linear isometry of $F$ into $H$.
\end{lemma}
With our new definition, this is restated in our project as follows:
\begin{lemma}\label{lem:32158}
    Let $E, F, H$ be normed spaces over $\bbK$, with $H$ spherically complete. If there exists a linear isometry $f\colon E\lto H$, then it factors through any immediate map $E\lto F$.
    $$\begin{tikzcd}
    F\ar[dr,dashed]&\\
    E\ar[u,"g"]\ar[r,"f"']&H
    \end{tikzcd}$$
\end{lemma}
\begin{leancode}[hyperurl={https://github.com/YijunYuan/SphericalCompleteness/blob/717759d2b22a5720e74c14160a8385ac17f32abc/SphericalCompleteness/NormedVectorSpace/Immediate.lean\#L162-L195}]
theorem IsImmediate.exists_linearIsometry_comp_eq {H : Type*} 
    [NormedAddCommGroup H] [NormedSpace 𝕜 H] [IsUltrametricDist H] 
    [SphericallyCompleteSpace H]
    (f : E →ₗᵢ[𝕜] F) (hf : IsImmediate f) (g : E →ₗᵢ[𝕜] H) :
    ∃ (h : F →ₗᵢ[𝕜] H), LinearIsometry.comp (h : F →ₗᵢ[𝕜] H) (f : E →ₗᵢ[𝕜] F) = g
\end{leancode}
\begin{remark}
    This lemma is a quick consequence of the Hahn–Banach theorem (\Cref{thm:14044}). We do not elaborate its formalized proof here, for brevity.
\end{remark}

\subsection{Maximal immediate extensions: a Zorn's lemma argument}
In \cite{rooijNonArchimedeanFunctionalAnalysis1978}, the spherical completion is constructed via the following result:
\begin{proposition}\label{prop:30007}
    Let $E$ be a normed space over $\bbK$, which embeds into some normed space $F$ via a linear isometry $f\colon E\to F$. Consider the set $S$ of submodules $M$ of $F$ containing $\opn{Im}(f)$ such that the inclusion map $\opn{Im}(f)\lto M$ is immediate. Then
    \begin{enumerate}
        \item $S$ has a maximal element $\breve{E}$ with respect to inclusion.
        \item if $F$ is spherically complete, then $\breve{E}$ is a spherical completion of $E$.
    \end{enumerate}
\end{proposition}
The set $S$ is formalized as follows:
\begin{leancode}[hyperurl={https://github.com/YijunYuan/SphericalCompleteness/blob/717759d2b22a5720e74c14160a8385ac17f32abc/SphericalCompleteness/NormedVectorSpace/Immediate.lean\#L214-L231}]
def IsImmediate.immediateExtensionSubmodules 
    (f : E →ₗᵢ[𝕜] F) : Set (Submodule 𝕜 F) := 
    {M : Submodule 𝕜 F | ∃ hc : f.range ≤ M,
      ∀ v : M, (v ⟂ₘ LinearMap.range (Submodule.inclusion hc)) → v = 0 }
\end{leancode}
After this definition, we apply Zorn's lemma to show the existence of maximal elements in \lean{IsImmediate.immediateExtensionSubmodules}. Given a chain $C$ in this set, an upper bound $M_{\max}$ is given by the union of all elements of $C$. We verify that $M_{\max}$ lies in this set by showing that the inclusion map $\opn{Im}(f)\lto M_{\max}$ is immediate. This yields the formalized proof of the first assertion of \Cref{prop:30007}:
\begin{leancode}[hyperurl={https://github.com/YijunYuan/SphericalCompleteness/blob/717759d2b22a5720e74c14160a8385ac17f32abc/SphericalCompleteness/NormedVectorSpace/Immediate.lean\#L266-L310}]
theorem IsImmediate.exists_maximal_immediateExtensionSubmodule (f : E →ₗᵢ[𝕜] F) :
    ∃ m, Maximal (fun x ↦ x ∈ immediateExtensionSubmodules E F f) m 
\end{leancode}
We define \lean{maximalImmediateExtensionEmbedding f} to be the inclusion map of $E$ into its maximal immediate extension in $F$, and prove that it gives a spherical completion of $E$ when $F$ is spherically complete:
\begin{leancode}[hyperurl={https://github.com/YijunYuan/SphericalCompleteness/blob/717759d2b22a5720e74c14160a8385ac17f32abc/SphericalCompleteness/NormedVectorSpace/SphericalCompletion/Basic.lean\#L61-L107}]
instance instIsSphericalCompletionOfMaximalImmediateExtensionSubmodule
    [SphericallyCompleteSpace F] (f : E →ₗᵢ[𝕜] F) :
    IsSphericalCompletion (maximalImmediateExtensionEmbedding f) 
\end{leancode}

\subsection{Existence of spherically complete extension}
With \Cref{prop:30007}, the existence of spherical completions of $E$ is reduced to the existence of a spherically complete extension of $E$. This is guaranteed by the following result:
\begin{proposition}[{cf. \cite[Theorem 4.1]{rooijNonArchimedeanFunctionalAnalysis1978}}]\label{prop:58585}
    Let $E$ be a normed space over $\bbK$. Let 
    $\ell^\infty(E)$ be the space of all bounded $\bfN$-indexed $E$-sequences, equipped with the supremum norm. Let
    $$c_0(E)\coloneqq \left\{(f_n)_{n\in\bfN}\in\ell^\infty(E)\colon \lvert f_n\rvert \to 0 \text{ as } n \to \infty\right\}.$$
    Then
    \begin{enumerate}
        \item $c_0(E)$ is a closed subspace of $\ell^\infty(E)$. In particular, $\ell^\infty(E)/c_0(E)$ is a normed space over $K$.
        \item $\ell^\infty(E)/c_0(E)$ is spherically complete.
        \item The diagonal embedding $i\colon E\to \ell^\infty(E)/c_0(E),\ x\mapsto [(x,x,x,\cdots)]$ is an isometric linear embedding.
    \end{enumerate}
    In particular, every normed space over $\bbK$ admits a spherically complete extension.
\end{proposition}
In Mathlib, $\ell^\infty(E)$ is already formalized as \lean{lp E ⊤}, which is a normed space over $K$. We formalize the subspace $c_0(E)$ as follows:
\begin{leancode}[hyperurl={https://github.com/YijunYuan/SphericalCompleteness/blob/717759d2b22a5720e74c14160a8385ac17f32abc/SphericalCompleteness/NormedVectorSpace/SphericalCompletion/LpQuotient.lean\#L58-L76}]
def c₀ (G : ℕ → Type*) [∀ i, NormedAddCommGroup (G i)]
    [∀ i, NormedSpace 𝕜 (G i)] : Submodule 𝕜 ↥(lp G ⊤) where
  carrier := {f : lp G ⊤ | Tendsto (fun n ↦ ‖f n‖) atTop (𝓝 0)} -- the underlying set
  add_mem' := ... -- closedness under addition. Omitted.
  zero_mem' := ... -- zero element is in c₀. Omitted.
  smul_mem' := ... -- closedness under scalar multiplication. Omitted.
\end{leancode}
The closedness of $c_0(E)$ and the resulting spherical completeness of the quotient space are not difficult to verify, and we implement them as two \lean{instance}s.

Finally, the diagonal embedding is formalized as a \lean{LinearIsometry}:
\begin{leancode}[hyperurl={https://github.com/YijunYuan/SphericalCompleteness/blob/717759d2b22a5720e74c14160a8385ac17f32abc/SphericalCompleteness/NormedVectorSpace/SphericalCompletion/LpQuotient.lean\#L291-L334}]
def lpQuotientEmbedding :
    E →ₗᵢ[𝕜] ((lp (fun (_ : ℕ) ↦ E) ⊤)⧸ c₀ 𝕜 (fun (_ : ℕ) ↦ E)) where
  toFun x := by
    refine (QuotientAddGroup.mk' _) (⟨fun (_ : ℕ) ↦ x, ?_⟩)
    ... -- Proof of the fact that the diagonal embedding belongs to $\ell^\infty$. Omitted.
  ... -- map_add', map_smul', norm_map'. Omitted.
\end{leancode}

\subsection{Weak universal properties of the spherical completion}
As we mentioned at the beginning of this section, the spherical completion does not enjoy as good a universal property as the usual completion of metric spaces. Still, some weaker universal-property-like results can be established.

To begin with, following van Rooij, we observe that the uniqueness of spherical completions in Fleischer's theorem (cf. \Cref{thm:63198}) follows from \Cref{lem:32158}: take two spherical completions $f_1\colon E\lto F_1$ and $f_2\colon E\lto F_2$ of $E$. Since $f_1$ is immediate and $F_2$ is spherically complete, \Cref{lem:32158} gives an isometric linear map $g\colon F_1\lto F_2$ such that $g\circ f_1=f_2$. As a result, the image of $g$ is a spherically complete subspace of $F_2$, which equals $F_2$ by the minimality of $F_2$. Hence $g$ is surjective, and consequently an isometric isomorphism.

On the other hand, the spherical completion provides another way to describe the spherical completeness of ultrametric normed spaces:
\begin{corollary}
    Let $E$ be a normed space over $\bbK$. Then the following are equivalent:
    \begin{enumerate}
        \item $E$ is spherically complete.
        \item The embedding of $E$ into its spherical completion $\breve{E}$ is surjective.
        \item $E$ is maximally complete, i.e. there is no non-trivial immediate extension of $E$.
    \end{enumerate}
\end{corollary}
\begin{leancode}[hyperurl={https://github.com/YijunYuan/SphericalCompleteness/blob/717759d2b22a5720e74c14160a8385ac17f32abc/SphericalCompleteness/NormedVectorSpace/SphericalCompletion/Basic.lean\#L233-L251}]
theorem IsSphericalCompletion.iff_embedding_to_sphericalCompletion_surjective
    {ι : E →ₗᵢ[𝕜] F} (hF : IsSphericalCompletion ι) :
    SphericallyCompleteSpace E ↔ Function.Surjective ι
\end{leancode}
\begin{leancode}[hyperurl={https://github.com/YijunYuan/SphericalCompleteness/blob/717759d2b22a5720e74c14160a8385ac17f32abc/SphericalCompleteness/NormedVectorSpace/SphericalCompletion/Basic.lean\#L257-L294}]
theorem iff_maximallyComplete :
    SphericallyCompleteSpace E ↔ MaximallyComplete 𝕜 E
\end{leancode}
The proof is not difficult, and we omit the details here.

Finally, we show that the spherical completion is minimal among all spherically complete extensions of $E$, which can be viewed as a weaker analogue of \Cref{thm:14845}:
\begin{corollary}\label{coro:23722}
    Let $E,F$ be normed spaces over $\bbK$. Let $f\colon E\lto F$ be a linear isometry and let $\iota\colon E\lto\breve{E}$ be a spherical completion of $E$. Then there exists a linear isometry $T\colon \breve{E}\lto F$ such that $T\circ \iota=f$.
\end{corollary}
\begin{leancode}[hyperurl={https://github.com/YijunYuan/SphericalCompleteness/blob/717759d2b22a5720e74c14160a8385ac17f32abc/SphericalCompleteness/NormedVectorSpace/SphericalCompletion/Basic.lean\#L218-L231}]
theorem IsSphericalCompletion.weak_universal_property
    {F' : Type*} [NormedAddCommGroup F'] [NormedSpace 𝕜 F'] [IsUltrametricDist F']
    [SphericallyCompleteSpace F'] (f : E →ₗᵢ[𝕜] F') 
    {ι : E →ₗᵢ[𝕜] F} (hF : IsSphericalCompletion ι) :
    ∃ T : F →ₗᵢ[𝕜] F', T.comp ι = f
\end{leancode}

\backmatter
\printbibliography
\end{document}